\documentclass[letterpaper,12pt]{article}

\usepackage{color}
\usepackage{amssymb}
\usepackage{amsmath}
\usepackage{amscd}
\usepackage{mathtools}
\usepackage[all]{xy}
\usepackage{tikz-cd}

\usepackage{hyperref}
\hypersetup{colorlinks=true,linkcolor=blue,anchorcolor=blue,citecolor=blue}
\usepackage{cleveref}

\usepackage{color}

\newtheorem{theo}{Theorem}[section]
\newtheorem{prop}[theo]{Proposition}
\newtheorem{lem}[theo]{Lemma}
\newtheorem{cor}[theo]{Corollary}
\newtheorem{defi}[theo]{Definition}
\newtheorem{rem}[theo]{Remark}
\newtheorem{conj}[theo]{Conjecture}

\def \nr {{\rm nr}}

\def \K{{\mathcal K}}
\def \deg {{\rm{deg}}}
\def \Br {{\rm{Br}}}
\def \Si {{\Sigma}}
\def \si {{\sigma}}
\def \Ga {{\Gamma}}

\def \Pic {{\rm {Pic}}}

\def \Gal {{\rm{Gal}}}
\def \Ker {{\rm{Ker}}}

\def \Frob{{\rm{Frob}}}
\def \AA{{\mathbf A}}
\def \A{{\mathbb A}}
\def \P{{\mathbb P}}
\def \Spec {{\rm{Spec}}}
\def \dim {{\rm{dim}}}
\def \Hom {{\rm {Hom}}}
\def \End {{\rm {End}}}

\def \PGL {{\rm {PGL}}}

\def \Aut{{\rm Aut}}

\def \SB{{\mathcal{S}}} 
\def \loc{{\rm loc}}

\def\ov{\overline}

\def \Z {{\mathbb Z}}
\def \Q {{\mathbb Q}}
\def \F {{\mathbb F}}

\def \AA {{\rm A}}

\def \DD {{\rm D}}

\def \FF {{\rm F}}

\def \Tr {{\rm{Tr}}}

\def \val {{\rm{val}}}

\def \Sel{{\rm Sel}}

\def\G{{\mathbb G}}

\def\GG{{\cal G}}

\def\T{{\cal T}}

\def\lra{\longrightarrow}

\def\N{{\rm N}}
\def\H{{\rm H}}

\def\Tr{{\rm Tr}}

\def\K{{\cal K}}

\def\inv{{\rm inv}}
\def\Kum{{\rm Kum}}

\def\NS{{\rm NS\,}}

\def\Si{\Sigma}
\def\si{\sigma}

\def\val{{\rm val}}
\def\Ram{{\rm Ram}}

\def\discr{{\rm discr\,}}
\def\Ga{\Gamma}

\def\e{\varepsilon}
\def\et{{\rm{\acute et}}}

\def\QQ{{\mathcal Q}}

\def\bA{{\bf A}}

\def\del{{\partial}}

\def\Kumm{{\mathcal{K}}}
\def\Pen{{\mathcal{P}}}

\DeclareFontEncoding{OT2}{}{} 
\DeclareTextFontCommand{\textcyr}{\fontencoding{OT2}\fontfamily{wncyr}\fontseries{m}\fontshape{n}\selectfont}

\newcommand{\Sha}{\textcyr{Sh}}

\newcommand{\bthe}{\begin{theo}}
\newcommand{\ble}{\begin{lem}}
\newcommand{\bpr}{\begin{prop}}
\newcommand{\bco}{\begin{cor}}
\newcommand{\bde}{\begin{defi}}
\newcommand{\brem}{\begin{rem}}
\newcommand{\bconj}{\begin{conj}}
\newcommand{\ethe}{\end{theo}}
\newcommand{\ele}{\end{lem}}
\newcommand{\epr}{\end{prop}}
\newcommand{\eco}{\end{cor}}
\newcommand{\ede}{\end{defi}}
\newcommand{\erem}{\end{rem}}
\newcommand{\econj}{\end{conj}}

\title{Hasse principle for intersections of two quadrics via Kummer surfaces}
\author{Adam Morgan and Alexei Skorobogatov}

\date{\today}
\begin{document}
\maketitle

\begin{abstract}
\noindent 
We prove new cases of the Hasse principle for Kummer surfaces constructed from 2-coverings
of Jacobians of genus 2 curves, assuming finiteness of relevant Tate--Shafarevich groups.
Under the same assumption, we deduce the Hasse principle for quartic del Pezzo surfaces with trivial 
Brauer group and irreducible or completely split characteristic polynomial, hence
the Hasse principle for smooth complete intersections of two quadrics in the projective space
of dimension at least 5.
\end{abstract}

{\small
\tableofcontents
}

\section*{Introduction}

The study of the Hasse principle for geometrically rational varieties has a long history.
Fundamental contributions to this subject were made by Manin who introduced the Brauer--Manin obstruction,
and by Colliot-Th\'el\`ene, Sansuc and Swinnerton-Dyer who pioneered the use of descent and
fibration methods.

By the Castelnuovo--Manin--Iskovskikh classification, smooth and projective
geometrically rational surfaces 
that are minimal in the sense of birational geometry (every birational morphism
onto a smooth surface must be an isomorphism) fall into two families: 
conic bundles and del Pezzo surfaces \cite[Theorem 1]{Isk79}.
The existence of a conic fibration greatly simplifies the study of rational points,
essentially because universal torsors over conic bundle surfaces have a simple
description in terms of certain complete intersections of quadrics.
As a result, under some assumptions on the degrees of singular closed fibres of the conic fibration,
one can prove that the Brauer--Manin
obstruction is the only obstruction to the Hasse principle and weak approximation
(see the survey \cite{CT92}, \cite[\S 4.4 and Chapter 7]{Sko01}, \cite{BMS}).
By contrast, no unconditional results of this kind are known for del Pezzo
surfaces without a conic fibration, if we exclude the well understood case of del
Pezzo surfaces of degree larger than $4$. 
The most interesting open cases are smooth cubic surfaces in $\P^3_k$
and quartic del Pezzo surfaces, which are smooth intersections of two quadrics in $\P^4_k$. 
(Universal torsors over arbitrary del Pezzo surfaces were described in \cite{SS07}, but
this description looks too difficult to work with directly.) 

Attempts were made to prove
the Hasse principle conditionally on the  conjecture of Tate and Shafarevich about finiteness
of the Tate--Shafarevich group of abelian varieties. 
This is known for many elliptic curves and for some higher-dimensional
abelian varieties. A conceptual reason comes from a theorem of Manin which says
that the Brauer--Manin obstruction is the only obstruction to the Hasse principle
on torsors for an abelian variety provided its Tate--Shafarevich group is finite, 
see \cite[Theorem 6.2.3]{Sko01}. Thus it seems reasonable
to build on this conjecture when dealing with surfaces that do not have a pencil of rational curves.

For {\em diagonal} cubic surfaces with coefficients 
satisfying certain divisibility properties, Swinnerton-Dyer proved
the Hasse principle conditionally on finiteness of the Tate--Shafarevich group of
elliptic curves \cite{SD01}. He crucially used the geometric fact that any diagonal cubic surface
is birational to a quotient of a torsor for an abelian surface, namely the product of two diagonal cubic curves.
Extending this approach, Swinnerton-Dyer's method was developed (and later much simplified)
to prove the Hasse principle for many
Kummer surfaces constructed from 2-coverings of abelian surfaces, assuming finiteness of the relevant
Tate--Shafarevich groups, see \cite{SSD05, HS16, Har19, Mor23b}. 

In his thesis, Wittenberg developed an earlier version of this method 
applicable to surfaces with a pencil of genus 1 curves \cite{CSS98a} to
prove the Hasse principle for many classes of quartic del Pezzo surfaces, including the most general case \cite[Th\'eor\`eme 3.2]{Wit07}. In addition to the Tate--Shafarevich conjecture 
he had to assume Schinzel's Hypothesis (H), a notoriously difficult conjecture on 
prime values of polynomials, deemed to be inaccessible in the current state of analytic number theory.

In this paper we return to the philosophy of \cite{SD01}. We study the Hasse principle
on quartic del Pezzo surfaces using their representation as quotients of Kummer surfaces
which are themselves quotients of torsors for the Jacobians of genus 2 curves. 

 For a quartic del Pezzo surface $X$
represented as the zero set of two quadratic forms in 5 variables $\Phi_1=\Phi_2=0$, one calls 
$\det(\Phi_1-t\Phi_2)\in k[t]$ the {\em characteristic polynomial} of $X$.
Our main result is the following theorem.

\medskip

\noindent{\bf Theorem A} {\em 
Let $k$ be a number field. Assume that for every genus $2$ curve
over every finite extension of $k$ the $2$-primary torsion subgroup of the 
Tate--Shafarevich group of its Jacobian is finite.
Then the quartic del Pezzo surfaces $X$ over $k$ whose characteristic polynomial is irreducible, or whose characteristic polynomial is
completely split and the natural map $\Br(k)\to\Br(X)$ is surjective, satisfy the Hasse principle.}

\medskip

To prove Theorem A, we use a new geometric idea that allows us to reduce the study
of rational points on 
quartic del Pezzo surfaces to Kummer surfaces attached to 2-coverings of Jacobians of genus 2 curves,
by exploiting a relation between these surfaces spelled out in \cite{S10}. 
Namely, we show that the projectivised tangent bundle
of a quartic del Pezzo surface $X$ is birationally equivalent to a family
of Kummer surfaces over $\P^1_k$, each of which is a double cover of $X$.
Although classical in spirit, this construction appears to be new. 
Applying the fibration
theorem for zero-cycles of Harpaz and Wittenberg \cite[Theorem 6.2]{HW16} to this family, we
prove the main unconditional result of this paper, Theorem \ref{dp}. 
It says that the Hasse principle with Brauer--Manin obstruction for $X$ over $k$ follows from
the Hasse principle with algebraic Brauer--Manin obstruction 
for Kummer surfaces defined over finite extensions of $k$ of odd degree and satisfying
certain local conditions. In the proof, one obtains a closed point of odd degree on $X$, but
since $X$ is an intersection of two quadrics, a theorem of Amer--Brumer
 implies that $X$ has a $k$-point. 
The required results on the existence of rational points on Kummer surfaces are proved
conditionally on finiteness of the Tate--Shafarevich groups of Jacobians of genus 2 curves. 
In the completely split case we use that the algebraic Brauer--Manin obstruction is 
the only obstruction to the Hasse principle on the relevant
Kummer surfaces, as proved by Harpaz in \cite[Theorem 1.3]{Har19}.
In the irreducible case, we establish the Hasse principle for the relevant Kummer surfaces in Theorem \ref{thm}, which we prove by extending the technique
of \cite{HS16, Har19, Mor23b} using work of the first named author \cite{Mor19, Mor23a}.

Although in this paper we only deal with quartic del Pezzo surfaces with
trivial Brauer group, our method is applicable to arbitrary quartic del Pezzo surfaces.
Instances where the Brauer group is non-trivial will be the subject of a forthcoming work of
Julian Lyczak and the first named author.

The following statement is a consequence of Theorem A, see \cite[\S 3.5.2]{Wit07} for
this implication.

\medskip

\noindent{\bf Theorem B} {\em
Let $k$ be a number field. Assume that $\Sha({\rm Jac}(C))[2^\infty]$ is finite
for every genus $2$ curve $C$ over every finite extension of $k$. Then the
smooth complete intersections of two quadrics in $\P^n_k$,  $n\geq 5$, satisfy the Hasse principle.
}

\medskip

Previously, this was known under the assumption of finiteness of the Tate--Shafarevich
groups of elliptic curves over number fields and Schinzel's Hypothesis (H), see
\cite[Th\'eor\`eme 3.3]{Wit07}. The best unconditional result in this direction is
due to Heath-Brown \cite{HB18}
who proved that smooth complete intersections of two quadrics in $\P^n_k$
satisfy the Hasse principle for $n\geq 7$. See \cite{CT24} for a different proof.

We briefly outline the structure of the paper. Section \ref{S1} contains various preliminary results.
In Section \ref{2} we introduce the main geometric construction of the paper, namely 
a stable birational equivalence between a quartic del Pezzo surface $X$
and a threefold $Z$ fibred into Kummer surfaces over $\P^1_k$. Using a natural
realisation of $Z/\P^1_k$ inside a fibration into 3-dimensional Severi--Brauer
varieties, we show the triviality of the vertical Brauer group of $Z/\P^1_k$. This result
is not used in the proof of Theorem A, but it is interesting in its own right.
In Section \ref{S3} we deduce the Hasse principle with Brauer--Manin obstruction for $X$ from
the Hasse principle with algebraic Brauer--Manin obstruction of a special type
for certain closed fibres of $Z/\P^1_k$, see Theorem \ref{dp} for a precise statement. 
This reduces Theorem A (=Theorem \ref{adrasan2}) to
the relevant cases of the Hasse principle for Kummer surfaces,
see Theorem \ref{thm}. Section \ref{S4} contains the proof of Theorem \ref{thm}.

The work on this paper started at the BIRS--CMI workshop ``New directions in rational points"
at the Chennai Mathematical Institute.
The authors are very grateful to Julian Lyczak for helpful discussions and for suggesting to use the relation between
quartic del Pezzo and Kummer surfaces to study rational points on quartic del Pezzo surfaces. 
We thank Olivier Wittenberg for explaining to us the relevant results of \cite{HW16}. 
The authors are extremely grateful to the referee for the thorough reading of the paper and 
very helpful suggestions.

Throughout this work, the first named author was supported by the Engineering and Physical Sciences Research Council (EPSRC) grant EP/V006541/1 ‘Selmer groups, Arithmetic Statistics and Parity Conjectures’.

\section{Preliminaries} \label{S1}

\subsection{Algebraic notation} \label{1.1}

Let $k$ be a field of characteristic zero
with algebraic closure $\bar k$ and absolute Galois group $\Ga=\Gal(\bar k/k)$. 

For a variety $X$ over $k$ we write $\Br(X)$ for the cohomological Brauer group 
$\H^2_\et(X,\G_m)$.
We use the standard notation $\Br_0(X)$ for the image of the natural map $\Br(k)\to\Br(X)$
and $\Br_1(X)$ for the kernel of the natural map $\Br(X)\to\Br(X_{\bar k})$.

\medskip

Let ${\rm S}_5$ be the group of permutations of $\{1,2,3,4,5\}$.
Define the ${\rm S}_5$-module $G$ as the zero-sum submodule of the permutation module 
$\mathbb{F}_2[{\rm S}_5/{\rm S}_4]\simeq (\F_2)^5$. Since $5$ is odd, we have a direct sum decomposition of ${\rm S}_5$-modules
$\mathbb{F}_2[{\rm S}_5/{\rm S}_4]\cong G \oplus \F_2$. 
It is clear that $G$ is a faithful ${\rm S}_5$-module.

In all of this paper, $P(t)$ is a monic separable polynomial in $k[t]$ of degree $5$. It gives rise
to an \'etale $k$-algebra $A=k[t]/(P(t))$ and an embedding of the $k$-scheme
$\Si=\Spec(A)$ into the affine line $\A^1_k$. We have $A=\oplus_{\si\in\Sigma} k_\si$, 
where $k_\si$ is the residue field of a closed point $\si\in\Sigma$.

Let $\Gal(P)$ be the Galois group of $P(t)$, 
and let $\varphi:\Ga\to \Gal(P)$ be the natural surjection. 
The action of $\Gal(P)$ on the roots of $P(t)$ gives rise to an embedding of 
$\Gal(P)$ into ${\rm S}_5$ (defined up to conjugation).
The Galois group $\Ga$ acts on $G$ via $\varphi$, making $G$ a finite commutative group 
$k$-scheme. Namely, we have 
$$G\cong R_{A/k}^1(\mu_2):=\Ker[R_{A/k}(\mu_2)\to\mu_2],$$
where $R_{A/k}$ is the Weil restriction of scalars and the map is induced by the norm map $A\to k$. 
Since $\dim_k(A)=5$ is odd, the exact sequence
$$1\to G\to R_{A/k}(\mu_2)\to\mu_2\to 1$$
is split, giving an isomorphism 
$$\H^1(k,G)\cong (A^\times/A^{\times 2})_1=\big(\prod_{\si\in\Sigma} 
k_\si^\times/k_\si^{\times 2}\big)_1,$$
where the subscript 1 denotes the subgroup of elements of norm $1\in k^\times/k^{\times 2}$.
Note that there is also an isomorphism $G\cong R_{A/k}(\mu_2)/\mu_2$, which induces
an isomorphism $\H^1(k,G)\cong A^\times/(k^\times A^{\times 2})$.

The fundamental datum in this paper is a pair $(P(t),\delta)$, where $\delta\in\H^1(k,G)$. 
We have $\delta=(\delta_\si)$, where
$\delta_\si\in k_\si^\times/k_\si^{\times2}$ are such that 
$$\prod_{\si\in\Sigma}{\rm N}_{k_\si/k}(\delta_\si)=1\in k^\times/k^{\times 2}.$$
We choose a lifting $\delta'\in A^\times$ of $\delta$ such that 
${\rm N}_{A/k}(\delta')\in k^{\times 2}$ and define
$\Delta:=\Spec(A[x]/(x^2-\delta')).$
It is clear that the isomorphism class of $\Delta$
does not depend on the lifting $\delta'\in A^\times$ of $\delta$.
The obvious morphism of $k$-schemes $\Delta\to\Sigma$
gives rise to a $\Ga$-equivariant 2-to-1 map $\Delta(\bar k)\to\Sigma(\bar k)$. We choose a set-theoretic section $s: \Sigma(\bar k) \to \Delta(\bar k)$ to this map. This identifies the group of permutations of the 10-element set $\Delta(\bar k)$ that respect 
the map $\Delta(\bar k)\to\Sigma(\bar k)$ with the semidirect product $(\Z/2)^5\rtimes {\rm S}_5$. The group $k$-scheme $R_{A/k}(\mu_2)$ acts faithfully on $\Delta$ preserving the fibres
of $\Delta\to\Sigma$, thus $R_{A/k}(\mu_2)(\bar k)$ is naturally identified with $(\Z/2)^5$ in
$(\Z/2)^5\rtimes {\rm S}_5$. It follows that $G(\bar k)=R_{A/k}^1(\mu_2)(\bar k)$ 
is identified with the zero sum subgroup $(\Z/2)^4\subset(\Z/2)^5$. 

For each $r\in \Sigma(\bar k)$, define a function $\tilde\delta_r:\Gamma \to \{\pm 1\}$ by setting $\tilde\delta_r(g)=1$ if $g s(r)=s(gr)$, and $\tilde\delta_r(g)=-1$ if $g s(r)\neq s(gr)$.

\ble \label{delta}
The function $\tilde \delta =(\tilde \delta_r):\Gamma \to (\Z/2)^5$ is a $1$-cocycle
$\Ga\to G(\bar k)$ with class~$\delta$. 
\ele
{\em Proof.} We have canonical isomorphisms of $\Ga$-sets 
$\Sigma(\bar k)=\Hom_{k-{\rm alg}}(A,\bar k)$ and 
$\Delta(\bar k)=\Hom_{k-{\rm alg}}(A',\bar k)$, where $A'=A[x]/(x^2-\delta')$.
Let $r\in \Sigma(\bar k)$.
The choice of a section $s: \Sigma(\bar k) \to \Delta(\bar k)$ is equivalent to the choice of a
lifting of each homomorphism of $k$-algebras $r\colon A\to\bar k$ 
 to a homomorphism of $k$-algebras $\tilde r\colon A'\to \bar k$, that is, to the choice of a 
square root of $r(\delta')$ in $\bar k$. By the definition of $\tilde\delta_r$ we have
\begin{equation}
\tilde\delta_r(g)=g(\sqrt{r(\delta')})/\sqrt{gr(\delta')}.\label{char}
\end{equation}
Since 
${\rm N}_{A/k}(\delta')\in k^{\times 2}$, we have $\prod_{r\in\Sigma(\bar k)}\tilde\delta_r=1$.
 For each closed point $\si\in\Sigma$, formula (\ref{char}) with $r\in \si$ describes
the character $\Gal(\bar k/k_\si)\to \{\pm 1\}$ of the extension
$k_\si(\sqrt{\delta'_\si})/k_\si$. Thus 
$\tilde\delta=(\tilde\delta_r)\colon\Ga\to G(\bar k)$ is a 1-cocycle with class $\delta=(\delta_\si)$. 
 \hfill $\Box$

\medskip

We fix a lifting $\delta'$, a set-theoretic section $s$ and hence the cocycle $\tilde \delta$, once and for all.

\medskip

Let us recall a useful construction from \cite[\S 3]{HS16}.
Suppose that we are given a finite subset $T=\{\alpha_1,...,\alpha_t\}\subset\H^1(k,G)$. 
Let us represent each $\alpha_i$ by a $1$-cocycle $\tilde\alpha_i:\Ga\to G(\bar k)$. 
Then the map $\varphi_T\colon\Ga\to G(\bar k)^T\rtimes \textup{Gal}(P)$ given by 
$$\varphi_T(g)=\big(\tilde\alpha_1(g),\ldots,\tilde\alpha_t(g),\varphi(g)\big)$$
is a homomorphism. 
Denote by $k_T$ the fixed field of $\Ker(\varphi_T)$, so that 
we have an injection $\textup{Gal}(k_T/k)\hookrightarrow G(\bar k)^T\rtimes \Gal(P)$. 
The extension $k_T/k$  contains the splitting field of $P(t)$ and is the smallest Galois extension of $k$ through which each of the cocycles $\tilde\alpha_1,\ldots,\tilde\alpha_t$ factors.  

From Lemma \ref{delta} we see that the natural action of $\Ga$ on $\Delta(\bar k)$ is via
$$\varphi_{\delta}\colon \Ga\hookrightarrow G(\bar k)\rtimes\Gal(P)
\hookrightarrow (\Z/2)^5\rtimes {\rm S}_5.$$

\ble   \label{delta zero}
The subgroup $\Gal(k_\delta/k)\subset(\Z/2)^4\rtimes{\rm S}_5$ is conjugate to a subgroup of 
${\rm S}_5$ in $(\Z/2)^5\rtimes{\rm S}_5$ if and only if $\delta=0$.
\ele
{\em Proof.} 
It is clear that $\delta=0$ if and only if the morphism $\Delta\to\Sigma$ has a section.
The stabiliser of a subset of $\Delta(\bar k)$ that maps bijectively onto $\Sigma(\bar k)$ is
$h{\rm S}_5h^{-1}$ for some $h\in(\Z/2)^5$. 
 \hfill $\Box$
 
\medskip
 
Let us denote by $\Lambda_\delta$ the $k$-torsor for $G$ defined by the cocycle $\tilde\delta$.
The isomorphism class of $\Lambda_\delta$ depends only on 
$\delta=[\Lambda_\delta]\in\H^1(k,G)$.
As above, the actions of $\Ga$ on $\Lambda_\delta(\bar k)$ and on $\Delta(\bar k)$
determine each other. In both cases, $\Ga$ acts through the faithfully acting finite group 
$\Gal(k_\delta/k)$. The relation between $\Lambda_\delta$ and $\Delta$ 
can also be described as follows: the image of $\delta$ in $\H^1(k,R_{A/k}(\mu_2))$ induced by
$G\hookrightarrow R_{A/k}(\mu_2)$ corresponds to $[\Delta]\in\H^1(A,\mu_2)$
under the isomorphism $\H^1(A,\mu_2)\cong \H^1(k,R_{A/k}(\mu_2))$ given by the Weil
restriction.

\medskip

Let $(\F_2)^{A,\delta}$ be the vector space over $\F_2$ with a basis whose vectors bijectively
correspond to the closed points $\si\in\Sigma$ such that 
$\delta_\si\neq 1\in k_\si^\times/k_\si^{\times 2}$. 
Let $\F_2\subset(\F_2)^{A,\delta}$ be the subspace spanned by the vector $(1,\ldots,1)$. 
Define ${\rm B}_{A,\delta}$ as the subgroup of $(\F_2)^{A,\delta}/\F_2$ represented by the 
vectors $\gamma=(\gamma_\si)\in (\F_2)^{A,\delta}$ such that 
$$\prod_{\si\in\Sigma}{\rm N}_{k_\si/k}(\delta_\si)^{\gamma_\si}=1\in k^\times/k^{\times 2}.$$
If $\delta=0$, we define ${\rm B}_{A,\delta}=0$.

\subsection{Euler identities}

Let $B=k[t]/(Q(t))$ where $Q(t)\in k[t]$ is a monic separable polynomial of degree 
$n\geq 1$. Let $\theta$ be the image of $t$ in $B$. Write
$Q(t)=\prod_{i=1}^n(t-\theta_i)$,
where $\theta_i\in \bar k$. A polynomial $R(t)\in k[t]$ of degree at most $n-1$ is uniquely 
determined by its values at $\theta_1,\ldots,\theta_n$, thus we have the Lagrange interpolation
formula
$$R(t)=\sum_{j=1}^nR(\theta_j)\frac{\prod_{i=1,\, i\neq j}(t-\theta_i)}{\prod_{i=1,\, i\neq j}(\theta_j-\theta_i)}
=\Tr_{B/k}\left(R(\theta)\frac{Q(t)}{Q'(\theta)(t-\theta)}\right),$$
where $Q'(t)$ is the derivative of $Q(t)$. Applying this to $R(t)=1,\ldots,t^{n-1}, t^n-Q(t)$ we obtain 
$$t^m=\Tr_{B/k}\left(\theta^m\frac{Q(t)}{Q'(\theta)(t-\theta)}\right), \ 
m\leq n-1, \quad
t^n-Q(t)=\Tr_{B/k}\left(\theta^n\frac{Q(t)}{Q'(\theta)(t-\theta)}\right).$$
Subtracting the expression for $t^m$ from the expression for $t^{m-1}$ multiplied by $t$,
we obtain Euler identities
\begin{equation}
\Tr_{B/k}\left(\frac{\theta^m}{Q'(\theta)}\right)=0, \ \ m=0,\ldots, n-2, \quad
\Tr_{B/k}\left(\frac{\theta^{n-1}}{Q'(\theta)}\right)=1, \label{Euler}
\end{equation}
see also \cite[Ch.~III, \S 6, Lemme 2]{CL}.

We shall need to consider quadratic forms
on the $k$-algebra $B$ given by
$$\Phi_\infty(u)=\Tr_{B/k}\left(\frac{u^2}{Q'(\theta)}\right), \quad
\Phi_b(u)=\Tr_{B/k}\left((b-\theta)\frac{u^2}{Q'(\theta)}\right),$$
where $b\in k$.

\ble \label{harp}
The determinants of the matrices of $\Phi_\infty$ and $\Phi_b$
in the basis $1,\ldots,\theta^{n-1}$ of the $k$-vector space $B$
are $(-1)^{n(n-1)/2}$ and $Q(b)(-1)^{n(n-1)/2}$, respectively.
\ele
{\em Proof.} The calculation for $\Phi_\infty$ follows from (\ref{Euler}), using that
its matrix has a triangular form, see \cite[{\em loc.~cit.}]{CL}.
The matrix of $\Phi_b$ is the product of the matrix of $\Phi_\infty$ and 
the matrix of the linear transformation
given by multiplication by $b-\theta\in B$.
The determinant of this matrix is the norm $\N_{B/k}(b-\theta)=Q(b)$. \hfill $\Box$

\subsection{Transitive subgroups of ${\rm S}_5$}

Recall that a transitive subgroup of ${\rm S}_5$ is isomorphic to the cyclic group
$\mathbb{Z}/5$, the dihedral group $\DD_{10}$, the Frobenius group $\FF_{20}$, the alternating group $\AA_5$, or the full symmetric group ${\rm S}_5$. Thus it
contains a $5$-cycle, and  contains a double transposition if it is not cyclic.
 
\ble \label{lemma:assumptions_a_and_b}
Let $\GG$ be a transitive subgroup of ${\rm S}_5$. Then 
$G$ is a simple and faithful $\GG$-module such that 
$\textup{End}_\GG(G)\cong\mathbb{F}_{2^r}$, where $r=4$ if $\GG\cong\Z/5$, 
$r=2$ if $\GG\cong \DD_{10}$, and $r=1$ otherwise.  
In all cases we have $\H^1(\GG,G)=0$.  
\ele
{\em Proof.} 
The polynomial $t^5-1\in\F_2[t]$ is the product of
irreducible factors $t-1$ and $t^4+t^3+t^2+t+1$, thus we have an isomorphism of rings
$\F_2[t]/(t^5-1)\cong\F_{16}\oplus\F_2$. This is also a direct sum of simple 
$\F_2[t]/(t^5-1)$-modules.

The ring $\F_2[t]/(t^5-1)$ is the group ring of the cyclic
group $\Z/5$ generated by a 5-cycle $\tau\in {\rm S}_5$. 
As a module over itself, it can be identified with $\F_2[{\rm S}_5/{\rm S}_4]$,
and then $G$ is identified with $\F_{16}$. In particular, $G$ is a simple $\Z/5$-module.

Taking $\tau\in\GG$ we conclude that $G$ is a simple $\GG$-module. 

The $\F_2$-algebra $\textup{End}_\GG(G)$ is the centraliser of $\GG$ in 
$\End_{\Z/5}(G)\cong\F_{16}$, so
it is the invariant subfield $(\F_{16})^\GG\cong\F_{2^r}$ where $r$ is $1$, $2$, or $4$. 
If $\GG\cong \mathbb{Z}/5$, then $r=4$. If $\GG\cong \DD_{10}$, then 
$(\F_{16})^\GG$ is the invariant subfield of a double transposition, so in this case $r=2$.
In all other cases, we have $r=1$.

Let ${\mathcal H}\subset \GG$ denote the stabiliser of an element of $\{1,2,3,4,5\}$. Then  we have 
$$\H^1(\GG,G)\oplus \textup{Hom}(\GG,\mathbb{F}_2)\cong 
\H^1(\GG,\F_2[{\rm S}_5/{\rm S}_4])\cong \textup{Hom}\big({\mathcal H},\F_2\big), $$
where the second isomorphism follows from Shapiro's lemma. Considering each transitive subgroup of 
${\rm S}_5$ in turn, we check from this that $\H^1(\GG,G)=0$. \hfill $\Box$

\ble \label{is_iso_vart}
Let $T$ be a finite subset of $\H^1(k,G)$. If $P(t)$ is irreducible, then $\varphi_T$ is 
surjective
 if and only if the elements of $T$ are $\mathbb{F}_q$-linearly independent, where 
$\F_q=\End_{\Gal(P)}(G)$.
When this is the case, every abelian subextension of $k_T/k$ is contained in the splitting field of $P(t)$.
\ele
\noindent {\em Proof.} 
By Lemma \ref{lemma:assumptions_a_and_b} we have $\H^1(\textup{Gal}(P),G)=0$. Given this, the first part is \cite[Proposition 3.6]{HS16} and the second part is  \cite[Corollary 3.9]{HS16}.
 \hfill $\Box$

\subsection{Severi--Brauer fibrations}

We shall need to consider the following situation. Let $\SB$ be a smooth, projective,
geometrically integral variety over $k$ with a dominant morphism
$\pi\colon \SB\to \P^1_k$ such that the generic fibre $\SB_{k(t)}$ is a Severi--Brauer variety.
Assume, moreover, that the restriction of $\pi$ to $\P^1_k\setminus\Sigma$ is smooth,
and that each closed fibre $\pi^{-1}(\si)$, $\si\in\Sigma$, consists
of two geometric components of multiplicity 1 individually defined over $k_\si(\sqrt{\e_\si})$
for some $\e_\si\in k_\si^\times$. Let 
$\e\in A^\times/A^{\times 2}=\prod_{\si\in\Sigma} k_\si^\times/k_\si^{\times 2}$
be the image of $(\e_\si)$.

We write $\Br_{\rm vert}(\SB)$ for the vertical
Brauer group associated to the morphism $\pi\colon \SB\to \P^1_k$.
The following statement is well-known;
we give the proofs for the convenience of the reader.

\bpr \label{due} Let $\alpha\in\Br(k(t))$ be the class of the Severi--Brauer variety $\SB_{k(t)}$.
Then we have the following statements.

{\rm (i)}
The residue of $\alpha$ at $\si\in\Sigma$ is the class of $\e_\si$ in $k_\si^\times/k_\si^{\times 2}$.
The residue of $\alpha$ at any closed point of $\P^1_k\setminus\Sigma$ is zero.
We have $\prod_{\si\in\Sigma}{\rm N}_{k_\si/k}(\e_\si)\in k^{\times 2}$.

{\rm (ii)} 
The natural map $\Br(k)\to\Br(\SB)$ is injective unless $\e= 1\in A^\times/A^{\times 2}$,
in which case its kernel is generated by $\alpha$.

{\rm (iii)} The inclusion $\Br_{\rm vert}(\SB)\subset\Br(\SB)$ is an equality.

{\rm (iv)} There is a natural isomorphism $\Br(\SB)/\Br_0(\SB)\cong {\rm B}_{A,\e}$.

{\rm (v)} There is a natural isomorphism $\Br(\SB)/\Br_0(\SB)\cong\H^1(k,\Pic(\SB_{\bar k}))$.
\epr
{\em Proof.} (i) 
By Lichtenbaum's theorem \cite[Theorem 5.4.11]{GSz}, the natural map
$\Br(k(t))\to\Br(\SB_{k(t)})$ is surjective with kernel generated by $\alpha$.
Thus $\alpha$ goes to zero under the restriction map to $\Br(k(\SB))$, hence the residue
of $\alpha$ at the generic point of every closed fibre $\SB_x$ of $\pi$ is trivial. 
If $x\notin\Sigma$, then $\SB_x$ is geometrically integral, hence the residue field $k(x)$ is separably
closed in $k(\SB_x)$. Using \cite[Lemma 11.1.2]{CTS21} we conclude that the residue of $\alpha$
at $x$ is zero. 
A similar argument gives that
the residue of $\alpha$ at a closed point $\si\in\Sigma$ belongs to the subgroup of 
$\H^1(k_\si,\Z/2)$ generated by the class of
$\e_\si$. The triviality of this residue in the case when
$\e_\si\notin k_\si^{\times 2}$ implies the existence of a Severi--Brauer
scheme over the local ring of $\si$ in $\P^1_k$ with generic fibre $\SB_{k(t)}$, but this
contradicts the birational invariance of $\e_\si$, see \cite[Corollary 2.3]{S15}. 
The product formula comes from the Faddeev reciprocity law \cite[Theorem 1.5.2]{CTS21}.

(ii) We adapt the argument in the proof of \cite[Lemma 11.3.3]{CTS21}.
Restricting to the generic point 
gives an injective map $\Br(\SB)\hookrightarrow\Br(k(\SB))$, since $S$ is smooth over $k$. 
The kernel of $\Br(k)\to\Br(\SB)$ is thus equal to the kernel of the composition
$$\Br(k)\hookrightarrow \Br(k(t))\to\Br(\SB_{k(t)})\hookrightarrow \Br(k(\SB)),$$
where the last map is injective since $\SB_{k(t)}$ is smooth over $k(t)$.
We know that the kernel of $\Br(k(t))\to\Br(\SB_{k(t)})$ is generated by $\alpha$.
By (i) and the Faddeev reciprocity law \cite[Theorem 1.5.2]{CTS21}, $\alpha$
belongs to $\Br(k)$ if and only if $\e=1\in A^\times/A^{\times2}$. 

(iii) This follows from the surjectivity of the natural map $\Br(k(t))\to\Br(\SB_{k(t)})$.

(iv) The notation ${\rm B}_{A,\e}$ makes sense by virtue of (i).
In view of Lichtenbaum's theorem, (iii) follows from (ii) and the standard calculation of 
the vertical Brauer group, see \cite[Proposition 11.1.8]{CTS21}. 

(v) The generic fibre of $\pi\colon \SB_{\bar k}\to\P^1_{\bar k}$ is a Severi--Brauer variety over 
$\bar k(t)$, so it is isomorphic to $\P^3_{\bar k(t)}$ by Tsen's theorem. 
Thus $\SB_{\bar k}$ is  rational, hence $\Br(\SB_{\bar k})=0$.  The claim
follows from \cite[Proposition 5.4.2]{CTS21} and \cite[Remark 5.4.3 (2)]{CTS21}.
\hfill $\Box$

\subsection{Number fields}

When $k$ is a number field, we denote by $O_k$ the ring of integers of $k$. 
For a set of places $S$ of $k$, let $O_{k,S}$ be the subring of $k$ consisting of the elements
that are integral ouside of $S$. We denote by
$k_v$ the completion of $k$ at a place $v$. When $v$ is non-archimedean, we denote
the ring of integers of $k_v$ by $O_v$ and the residue field of $O_v$ by $\F_v$.
Let us denote by
$$\loc_v\colon\H^1(k,G)\to\H^1(k_v,G)$$
the natural restriction map. Define
$$\H^1_\nr(k_v,G):=\Ker[\H^1(k_v,G)\to \H^1(k_v^\nr,G)],$$
where $k_v^\nr$ is the maximal unramified extension of $k_v$.

For a finite field extension $k'/k$ we write $\Ram(k'/k)$ for the set of primes of $k$ ramified in $k'$.

In all of this paper we write $S_0$ for a finite set of places of $k$ that contains all the `bad' 
places for the pair $(P(t),\delta)$:
the archimedean places, the places above the prime 2, the places $v$
where $P(t)\not\in O_v[t]$, and the places dividing $\discr(P(t))$. We arrange that $\delta$
is unramified away from $S_0$, so that for $v\notin S_0$ we have
$\loc_v(\delta)\in \H^1_\nr(k_v,G)$.

For a place $v\notin S_0$ we have a Frobenius element ${\rm Frob}_v\in \Gal(k_\delta/k)$
defined up to conjugation. 

Increasing $S_0$ by a fixed finite set of primes, we ensure that
any quartic del Pezzo surface over $\F_v$ has more $\F_v$-points than any genus 1 curve.

\ble \label{lem:existence_of_extensions}
Assume that $P(t)$ is irreducible and $\textup{End}_{\textup{Gal}(P)}(G)=\mathbb{F}_q$.
Let $T$ be a finite $\mathbb{F}_q$-linearly independent subset of $\H^1(k,G)$.
Let $S$ be a finite set of places of $k$ containing $S_0\cup\Ram(k_T/k)$. Let $\sigma\in \textup{Gal}(k_T/k)$
be such that every homomorphism $\textup{Gal}(P)\to \mathbb{Z}/2$ vanishes on the image of 
$\sigma$.  Then there are a prime $\mathfrak{p}\notin S$ such that 
$\Frob_{\mathfrak{p}}\in \Gal(k_T/k)$ lies in the conjugacy class of $\sigma$, and a  quadratic extension $F/k$ ramified at $\mathfrak{p}$, split at all places of $S$ and unramified at all other primes. 
 \ele
 \noindent {\em Proof.} 
 Let $\mathfrak{m}$ be the formal product of $8$ and all places in $S$. Let $k_\mathfrak{m}$ be the ray class field associated to the modulus $\mathfrak{m}$, and let $k_{\mathfrak{m},2}/k$ denote its maximal multiquadratic subextension.  By Lemma \ref{is_iso_vart}, the intersection of $k_{\mathfrak{m},2}/k$ with $k_T/k$ is contained in the splitting field of $P(t)$. Let $L$ be the compositum of $k_T$ and $k_{\mathfrak{m},2}$. Since every homomorphism $\textup{Gal}(P)\to \mathbb{Z}/2$ vanishes on the image of $\sigma$, we can find $\sigma'\in \textup{Gal}(L/k)$ with image $\sigma$ in $\textup{Gal}(k_T/k)$ and trivial image in $\textup{Gal}(k_{\mathfrak{m},2}/k)$. 
 
 Let $\mathfrak{p}$ be any prime outside $S$ such that $\textup{Frob}_{\mathfrak{p}}\in \textup{Gal}(L/k)$ lies in the conjugacy class of $\sigma'$. By construction, the image of $\textup{Frob}_{\mathfrak{p}}$ in $\textup{Gal}(k_T/k)$ lies in the conjugacy class of $\sigma$. Since the image of $\textup{Frob}_{\mathfrak{p}}$ in $\textup{Gal}(k_{\mathfrak{m},2}/k)$ is trivial, there is a fractional ideal $I$ supported on primes outside $S$ such that 
$\mathfrak{p}I^2=(a)$ for some totally positive element $a \in k^\times$  congruent to $1$ modulo
$8$ and modulo all odd primes in $S$. 
Using Hensel's lemma for the primes above 2, we check that
$F=k(\sqrt{a})$ satisfies the requirements of the statement.
 \hfill $\Box$

\subsection{Fibration theorem for zero-cycles}

The following proposition is a particular case of the general fibration theorem
for zero-cycles of Harpaz and Wittenberg
\cite[Theorem 6.2]{HW16}. The statement is adapted for our needs in this paper
and does not aim at maximum generality.

\bpr \label{h-w}
Let $k$ be a number field.
Let $Z$ be a smooth, projective, geometrically integral variety.
Let $f\colon Z\to \P^1_k$ be a dominant morphism with geometrically integral generic fibre
such that every geometric fibre of $f$ contains an irreducible component
of multiplicity $1$. Let $U\subset\P^1_k$ be a non-empty open subset such that
the restriction of $f$ to $Z_U=f^{-1}(U)$ is a smooth morphism $Z_U\to U$.
Let $B\subset \Br(Z_U)$ be a finite subgroup and let $m$ be any positive integer. 
Assume that there is an adelic point $(x_v)\in Z(\bA_k)$
orthogonal to $(B+f^*\Br(k(\P^1_k)))\cap \Br(Z)$. 
Then $U$ contains a closed point $R=\Spec(K)$, of degree coprime to $m$, such that $Z_R(\bA_K)^B\neq\emptyset$.
Moreover, for any non-archimedean places $v_1,...,v_n$ of $k$, we can choose $R$  such that, for each $i=1,...,n$, there is a place $w_i$ of $K$ over $v_i$ with $K_{w_i} \cong k_{v_i}$ and $R$ 
is arbitrarily close to $f(x_{v_i})$ in $\A^1_k(K_{w_i})$. 
\epr
{\em Proof.} Note that $(B+f^*\Br(k(\P^1_k)))\cap \Br(Z)$ is finite modulo $\Br(k)$, 
by \cite[Lemma 5.3]{HW16},
hence  is generated by finitely many elements $\alpha_1,\ldots,\alpha_r$ modulo $\Br(k)$. 
Each of $\alpha_1,\ldots,\alpha_r$ has trivial evaluation map at $Z(k_v)$ for almost all places $v$. 

Choose an arbitrary closed point $N$ in $Z_U$.

For each real place $v$, by a version of Bertini's theorem we can choose a smooth, projective,
geometrically irreducible curve $C_v\subset Z_{k_v}$ passing
through $x_v$ and $N$ such that $f(C_v)=\P^1_{k_v}$.
We choose a positive integer $r$ large enough so that if $y\in\Pic(\P^1_k)$ is a divisor
class of degree $1+m r\, \deg_k(N)$, then 
the inequality in condition (ii) of \cite[Theorem 5.1]{HW16} holds 
and the inequality in condition (iii) of \cite[Theorem 5.1]{HW16}
holds for every real place $v$ of $k$. 

For each place $v$, define $z_v$ to be an effective zero-cycle in $Z_{k_v}$
that is a small deformation of $x_v+mrN$ such that $f(z_v)$ is reduced and is contained in $U_{k_v}$.
For each real place $v$ we take $z_v$ to lie on $C_v$.
Since $\alpha_1,\ldots,\alpha_r$ have trivial evaluation maps at $Z(k_v)$ for almost all places $v$,  
at these places we can take $z_v$ to be
an {\em arbitrary} effective zero-cycle of degree  $1+mr\,\deg_k(N)$
such that $f(z_v)$ is reduced and contained in $U_{k_v}$,
and so arrange that $z_v$
is integral with respect to some integral model of $Z_U$. 
(An integral 
model of $Z_U$ is a faithfully flat $O_{k,S}$-scheme of finite type 
with generic fibre $Z_U$, where $S$ is a finite set of primes of $k$. A zero-cycle $z_v$ is integral 
if there is an integral model of $Z_U$ such that
each closed point in the support of $z_v$ extends to a finite $O_v$-subscheme of the model.)

Thus, in the notation of \cite{HW16}, we have 
$(z_v)\in{\rm Z}_{0,{\bf A}}^{{\rm eff, red,} y}(Z_U)$,
where $y\in\Pic(\P^1_k)$ has degree $1+mr\,\deg_k(N)$. The adelic zero-cycle $(z_v)$
satisfies conditions (i), (ii), (iii) of \cite[Theorem 5.1]{HW16}. 
An application of
\cite[Theorems 5.1, 6.2]{HW16} gives a closed point $R$ as in the theorem. \hfill $\Box$

\section{Projective geometry of pencils of quadrics in $\P^4_k$} \label{2}

\subsection{Quartic del Pezzo surfaces} \label{2.1} \label{new geom}

Let $k$ be a field of characteristic zero.

A quartic del Pezzo surface $X$ over $k$ is a smooth complete intersection of 
two quadrics in $\P^4_k$.
We have $X=Q_1\cap Q_2$, where $Q_1$ and $Q_2$ are any two distinct
quadrics containing $X$. Recall that a pencil of quadrics in $\P^4_k$ is a projective
line in the projective space of all quadrics in $\P^4_k$.
We denote the pencil spanned by $Q_1$ and $Q_2$ by $\Pen$. 
Let $\Phi_1$ and $\Phi_2$ be quadratic forms in 5 variables defining $Q_1$ and $Q_2$,
respectively.

Let $\Sigma\subset\P^1_k$ be the closed $k$-subscheme of zeros of
$\det(\mu \Phi_1-\nu \Phi_2)\in k[\mu,\nu]$. The smoothness of $X$ implies that 
$A:=\H^0(\Sigma, O_\Sigma)$ is an \'etale $k$-algebra of degree $\dim_k(A)=5$.
Choose a $k$-point in $\P^1_k\setminus\Sigma$ and let 
$\A^1_k=\Spec(k[t])\subset\P^1_k$ be its complement.
There is a unique separable monic polynomial $P(t)$ such that $A=k[t]/(P(t))$. 
We have $A=\oplus_{\si\in\Sigma} k_\si$, where $k_\si$ is the residue field of 
a closed point $\si\in\Sigma$. Let $\theta_\si$ be the image of $t$ in $k_\si$, 
and let $\theta=(\theta_\si)$ be the image of $t$ in $A=\oplus_{\si\in\Sigma} k_\si$. 

Let $\QQ\subset \P^4_k\times_k\P^1_k$ be the closed subvariety given by
$\mu \Phi_1-\nu \Phi_2=0$, where $(\mu:\nu)$ are homogeneous coordinates in $\P^1_k$.
The projection $p\colon\QQ\to\P^1_k$ is a quadric bundle of relative dimension 3 
which is smooth over  $\P^1_k\setminus\Sigma$. The fibre $\QQ_m$, where
$m=(\mu:\nu)\in \P^1_k(\bar k)$, is the zero set of the quadratic form $\mu \Phi_1-\nu \Phi_2=0$.
The closed fibres $\QQ_\si$, for $\si\in\Sigma$, are given
by quadratic forms of rank 4. The determinants of these quadratic forms are well-defined
up to a square, thus we have well-defined elements 
$\delta_\si\in k_\si^\times/k_\si^{\times2}$. Let $\delta=(\delta_\si)\in A^\times/A^{\times2}$.
We can lift $\delta$ to an element $\delta'\in A^\times=\prod_{\si\in\Sigma}k_\si^\times$.
It can be shown that $\delta\in (A^\times/A^{\times 2})_1$, so that
$N_{A/k}(\delta')\in k^{\times 2}$, see Lemma~\ref{duedue} below.
Thus we have $\delta\in\H^1(k,G)$. 

By \cite[Theorem 2.3]{S10} associating a pair $(\Sigma,\delta)$ to $X$
gives a bijection between isomorphism classes of quartic del Pezzo surfaces over $k$
and pairs $(\Sigma,\delta)$, where $\Sigma$ is a reduced 0-dimensional closed
subscheme of $\P^1_k$ of degree 5 considered up to the action of $\PGL(2,k)$, and
$\delta\in\H^1(k,G)$. 

Given a pair $(\Sigma,\delta)$, we can recover $X$ as the zero set of certain canonical equations.
Choose a coordinate $u$ on $A$;
it is well-defined up to multiplication by an element of $A^\times$. By \cite[\S 2]{S10}
the del Pezzo surface $X$ is isomorphic to the closed subscheme of $\P(R_{A/k}(\A^1_A))$ 
defined by the equations
\begin{equation}
\Tr_{A/k}\left(\frac{1}{P'(\theta)}\delta' u^2\right)=\Tr_{A/k}\left(\frac{\theta}{P'(\theta)}\delta' u^2\right)=0.
\label{dp4}
\end{equation}
Let $\Phi_\infty(u)$ be the first quadratic form in (\ref{dp4}). Its zero set is the quadric
$\QQ_\infty\subset\P^4_k$, which is smooth by Lemma \ref{harp}.
Every $\bar k$-quadric in the pencil $\Pen$ other than $\QQ_\infty$ is
the quadric $\QQ_b\subset\P^4_k$ given by the quadratic form
$$\Phi_b(u):=\Tr_{A/k}\left(\frac{b-\theta}{P'(\theta)}\delta' u^2\right)$$
for some $b\in \bar k$. 
By Lemma \ref{harp}, the quadric $\QQ_b$ is smooth if and only if $P(b)\neq 0$.

\medskip

Let us call $X_0\subset\P^4_k$ the surface given by (\ref{dp4}) with $\delta'=1$.
In fact, $X_0$ is isomorphic to the blowing-up of the image of 
$\Sigma\subset\P^1_k$ under the Veronese embedding $\P^1_k\hookrightarrow\P^2_k$. 
The strict transform of the image of $\P^1_k$ is a $k$-line $\lambda_0\subset X_0$ with 
homogeneous coordinates $(r:s)$ given by $u=r+s\theta$, cf.~\cite[Remark, p.~76]{S10}.

Equations (\ref{dp4}) show that the action of $R_{A/k}(\G_m)$ on $R_{A/k}(\A^1_A)$ 
by multiplication gives rise
to an action of $G\cong R_{A/k}(\mu_2)/\mu_2$ on $X$ and $X_0$ by automorphisms.
The surface $X$ is the twist of $X_0$ by the 1-cocycle $\tilde\delta:\Ga\to G(\bar k)$.
The action of $G$ on the $k$-scheme of 16 lines of $X_{\bar k}$
makes it a $k$-torsor for $G$ isomorphic to $\Lambda_\delta$. (Recall that $\delta=
[\Lambda_\delta]\in\H^1(k,G)$.)
The action of $G$ on the $k$-scheme of 10 pencils of conics of $X_{\bar k}$ identifies it with the
$G$-scheme $\Delta$. (See the end of Section \ref{one} for more details.)

The following well-known statement is given here for the ease of reference.

\bpr \label{coh dp}
Let $X$ be a quartic del Pezzo surface defined by a pair $(P(t),\delta)$. Then
$\H^1(k,\Pic(X_{\bar k}))\cong {\rm B}_{A,\delta}$.
\epr
{\em Proof.} Let $K=k(X)$. The generic point of $X$ gives rise to a $K$-point on $X_K=X\times_kK$
not contained in any of the 16 lines. Projecting from this point gives a birational
equivalence of $X_K$ with a cubic surface $Y$ over $K$ containing a $K$-line
$\ell\cong\P^1_K$. Residual conics of the plane sections through $\ell$ make $Y$
a conic bundle over $\P^1_K$ satisfying the assumptions of Proposition \ref{due} with $\e=\delta$.
Thus $\H^1(K,\Pic(Y_{\ov K}))\cong{\rm B}_{A,\delta}$. 
The natural isomorphism of $\Gal(\ov K/K)$-modules $\Pic(Y_{\ov K})\cong \Pic(X_{\ov K})\oplus\Z$
implies that $\H^1(K,\Pic(X_{\ov K}))\cong{\rm B}_{A,\delta}$. 
Using that $k$ is separably closed in $K$, we deduce the required isomorphism. \hfill $\Box$

\subsection{The associated Kummer fibration} \label{one}

Let $\P(\T_X)$ be the projectivisation of the tangent bundle on $X$. Thus $\P(\T_X)$ has
a smooth projective morphism to $X$ all of whose fibres are projective lines. In particular, $\P(\T_X)$
is birationally equivalent to $X\times_k\P^1_k$.

Let us write ${\rm Gr}(\P^1,\P^4)$ for the Grassmannian of projective lines in $\P^4_k$.
Define $Y\subset {\rm Gr}(\P^1,\P^4)\times\Pen$ to be the variety parameterising pairs 
$(\ell,m)$, where $\ell\subset\P^4_{\bar k}$ is a tangent line to $X$
contained in the quadric $\QQ_m$, where $m\in\Pen(\bar k)$. 

Let $Z\subset \P(\T_X)\times\Pen$ be the variety of triples $(x,\ell,m)$, where $x\in X(\bar k)$,
$\ell\subset\P^4_{\bar k}$ is a line tangent to $X$ at $x$
 contained in the quadric $\QQ_m$, where $m\in\Pen(\bar k)$. 
Note that the restriction of $\P(\T_X)\to X$ to a $\bar k$-line $\lambda\subset X$ 
has a canonical section sending a point $x\in \lambda(\bar k)$ to $(x,\lambda)$. 
This gives rise to 16 {\em disjoint} $\bar k$-lines 
$\tilde\lambda$ in $\P(\T_X)$ which project onto the 16 $\bar k$-lines $\lambda$ in $X$.
The inverse image of each $\tilde\lambda$ in $Z$ is $\tilde\lambda\times\Pen$, because
$\lambda$ is contained in every quadric of the pencil $\Pen$.

When $\ell$ is a tangent line to $X$ not contained in $X$,
there is exactly one quadric in $\Pen$ through $\ell$, thus 
the restriction of the projection $Z\to \P(\T_X)$ to 
$\P(\T_X)\setminus \coprod\tilde\lambda$ is an isomorphism. 
Associating to a point $(x,\ell)$ of $\P(\T_X)\setminus\coprod\tilde\lambda$
the unique quadric in $\Pen$ that contains $\ell$
defines a surjective morphism $\P(\T_X)\setminus\coprod\tilde\lambda\to \P^1_k$ that extends to 
a dominant rational map from $\P(\T_X)$ to $\P^1_k$. From the definition of $Z$
we see that $Z$ resolves the indeterminacy of this rational map, 
so $Z$ is isomorphic to the blowing-up 
of $\P(\T_X)$ in the closed subvariety 
$\coprod\tilde\lambda$. Since $\P(\T_X)$ and $\coprod\tilde\lambda$
are smooth, we conclude that $Z$ is smooth.

The projection $f\colon Z\to\Pen\cong\P^1_k$
sending $(x,\ell,m)$ to $m$ factors through
the projection $Y\to \Pen\cong\P^1_k$. The morphism $f$
restricted to $\tilde\lambda\times\Pen$ is the projection to $\Pen\cong\P^1_k$. We
thus have the following commutative diagram:
$$\xymatrix{(\coprod\tilde\lambda) \times\Pen\ar@{^{(}->}[r]\ar[d]&
Z\ar[r]\ar[d]&Y\ar[r]\ar[d]&\Pen\cong\P^1_k\\
\coprod\tilde\lambda \ar@{^{(}->}[r]\ar[d]&\P(\T_X)\ar[r]\ar[d]&{\rm Gr}(\P^1,\P^4)& \\
\cup\lambda \ar@{^{(}->}[r]&X&&}$$
Note that the action of $G$ preserves all varieties in this diagram and all morphisms
are $G$-equivariant.

\medskip

Let $U=\A^1_k\setminus\Sigma$. 
Let $C\to U$ be a smooth projective morphism
whose fibres are genus 2 curves with affine equation
$y^2=(x-t)P(x)$, where $t$ is the coordinate on $\A^1_k$ fixed above.
Let $J\to U$ be the relative Jacobian of $C\to U$.
A crucial property for us is that
the 2-torsion group subscheme $J[2]$ is constant: it is the pullback to
$U$ of the group $k$-scheme $G$. 

A $k$-scheme on which $G$ acts by automorphisms can be twisted by the 1-cocycle 
$\tilde\delta\colon \Gal(\bar k/k)\to G(\bar k)$.
Thus we obtain a smooth projective morphism $J_\delta\to U$ whose fibres are 2-coverings
of the fibres of $J\to U$. Passing to the fibrewise quotient by the antipodal involution we get
a projective morphism ${\mathcal K}'_\delta\to U$
whose fibres are singular Kummer surfaces with 16 singular points.
Its fibrewise desingularisation is a smooth projective morphism $\K_\delta\to U$
whose fibres are desingularised Kummer surfaces with 16 lines. 
The scheme parameterising the 16 lines in the fibres is $\Lambda_\delta\times_k U$.
When $\delta\in \H^1(k,G)$ is zero, we denote $\K_\delta$ by $\K$.

Let us write $Z_U=Z\times_{\P^1_k}U$.
The following theorem gives the main geometric construction on which this paper is based.

\bthe \label{t1}
There is an isomorphism of $U$-schemes $Z_U\cong \K_\delta$.
In particular, $\K_\delta$ is birationally equivalent to $X\times_k\P^1_k$.
\ethe
{\em Proof.} Recall that $t$ is a coordinate on $U\subset \A^1_k$.
It is easy to deduce from \cite[Theorem 3.1, Corollary 3.3]{S10} that $\K_\delta$
is the closed subset of 
$\P(\A^1_k\times R_{A/k}(\A^1_A))\times U$
given by equations (\ref{dp4}) together with the equation
\begin{equation}
v^2=\Tr_{A/k}\left(\frac{P(t)}{(t-\theta)P'(\theta)}\delta' u^2\right), \label{double}
\end{equation}
where $u$ is a coordinate in $R_{A/k}(\A^1_A)$ and $v$ is a new variable.
Forgetting the coordinate $v$ gives rise to a double cover $\K_\delta\to X\times U$, 
cf.~\cite[Proposition 3.2]{S10}. This double cover is $G$-equivariant, 
sends the 16 lines in each fibre of  $\K_\delta\to U$ isomorphically onto the 16 lines on $X$,
and is the twist by $\tilde\delta$ of the double cover $\K\to X_0\times U$.

We now show that $Z_U\to X\times_kU$ can be given by the same equation (\ref{double})
as $\K_\delta$. The following short argument suggested by the referee replaces a more
roundabout argument from the first version of this paper. Without loss of generality we can assume 
$\delta'=1$; the general case follows from this by twisting by the 1-cocycle $\tilde \delta$.

We note that $Z_U$ is the variety of triples $(x,\ell,t)$, where 
$x\in X_0(\bar k)$, $t\in\bar k$, $P(t)\neq 0$, and 
$\ell\subset\QQ_t$ is a line tangent to the quadric $\QQ_\infty$ at $x$. 
Let $u$ be a $\bar k$-point of $R_{A/k}(\A^1_A)$ above $x$.
The tangent hyperplane to $\QQ_\infty$ at $x$ is given by $\Tr_{A/k}(uw/P'(\theta))=0$.
Let $W$ be the hyperplane in $R_{A/k}(\A^1_A)$ given by this equation.
Thus $W$ is the orthogonal complement to $u$ with respect to the bilinear form
of $\Phi_\infty(u)=\Tr_{A/k}(u^2/P'(\theta))$. Equivalently, $W$ is the orthogonal complement to 
$u/(t-\theta)$ with respect to the bilinear form of $\Phi_t(u)=\Tr_{A/k}((t-\theta)u^2/P'(\theta))$. 
The fibre of $Z_U\to X_0\times_k U$ at $(x,t)$ consists of the two $\bar k$-lines in the 
2-dimensional quadric $\QQ_t\cap\P(W)$ through $x$, one from each of the two families of lines.
Thus the double cover $Z_U\to X_0\times_k U$ is given by extracting the square root
of $\discr(\QQ_t|_W)$.
By Lemma \ref{harp}, we have $\discr(\QQ_t)=P(t)$.
Up to multiplication by a square, the discriminant of the restriction of $\QQ_t$ 
to the orthogonal complement to $u/(t-\theta)$ is 
$$P(t)\Phi_t\left(\frac{u}{t-\theta}\right)=\Tr_{A/k}\left(\frac{P(t)}{(t-\theta)P'(\theta)}u^2\right).$$
This is precisely the right hand side of (\ref{double}) for $\delta'=1$. \hfill $\Box$

\brem{\rm
We have chosen $\A^1_k$ in the beginning of Section \ref{2.1}
as the complement to an arbitrary $k$-point in $\P^1_k\setminus\Sigma$. Applying
Theorem \ref{t1} with a different choice of $\A^1_k$
shows that the fibre of $f\colon Z\to\P^1_k$ at infinity is a smooth Kummer surface,
thus the restriction of $f\colon Z\to\P^1_k$ to $\P^1_k\setminus\Sigma$ is smooth.}
\erem

We now analyse the structure of the fibre $Z_\si$, where $\si$ is a closed point of $\Sigma$. 
The quadric $\QQ_\si$ is given by a quadratic form of rank 4 of discriminant 
$\delta_\si\in k_{k(\si)}^\times$ (well-defined up to multiplication  by a square), so it is a cone
over a non-degenerate quadric $\QQ_\si^\sharp\subset\P^3_{k(\si)}$.
Let $v_\si$ be the vertex of $\QQ_\si$, and let $\pi$ be the projection 
$\QQ_\si\setminus\{v_\si\}\to\QQ_\si^\sharp$ from $v_\si$.
Let $\phi^+$ and $\phi^-$ be the two families of $\bar k$-lines on $\QQ_\si^\sharp$.
Each of $\phi^+$ and $\phi^-$ is defined over $k_\si(\sqrt{\delta_\si})$.
As we vary $\si$, these families give rise to the ten families of conics on $X$ which are parameterised
by the $\bar k$-points of $\Delta$.  

The quadric $\QQ_\si$ has two families of lines, $F^+$ and $F^-$, consisting
of the lines whose image under $\pi$ belongs to a line in $\phi^+$ or in $\phi^-$, respectively.
Each of $F^+$ and $F^-$ is defined over $k_\si(\sqrt{\delta_\si})$.
The intersection $F^+\cap F^-$ consists of the lines through $v_\si$.

The morphism $\pi\colon X\to \QQ_\si^\sharp$ is a double cover ramified in 
a smooth curve $E_\si$ of genus 1 (a hyperplane section of $X$). 
If $x$ is a point of $X\setminus E_\si$, then there is 
exactly one line in $F^+$ which is tangent to $X$ at $x$ (respectively, in $F^-$).

Over $k_\si(\sqrt{\delta_\si})$, we have two projections $\QQ_\si^\sharp\to \P^1$.
Each of them represents $E_\si$ as a double cover of $\P^1$ ramified in four points.
Each of these four points is the intersection point of a pair of lines contained in $X$.
When $x$ is a point of $E_\si$ other than these eight points (four for each projection),
then the line $(v_\si x)$ is the only line in $\QQ_\si$ tangent to $X$ at $x$.

Finally, when $x$ is one of the four points, say defined by the family $\phi^+$, 
the tangent plane $T_{X,x}$ passes through $v_\si$ 
and projects onto a line of $\phi^+$. In this case, the family of lines in $F^+$ tangent to
$X$ at $x$ is parameterised by the projective line $\P(T_{X,x})$.

We conclude that the irreducible components of the fibres
of $f\colon Z_{\bar k}\to\P^1_{\bar k}$ over the points of $\Sigma(\bar k)$ correspond bijectively
and Galois-equivariantly  to the points of $\Delta(\bar k)$.
 Each of these is isomorphic to $X_{\bar k}$ blown-up
at four points, namely the singular points of the four singular conics in the family of conics 
given by a point of  $\Delta(\bar k)$. The two components of the fibre above each point
of $\Sigma(\bar k)$ intersect in a smooth curve of genus 1.

\subsection{The associated Severi--Brauer fibration} \label{two}

Recall that for a non-degenerate quadratic form $\Phi$ in 5 variables, the even Clifford algebra
$C^+(\Phi)$ is a central simple algebra over $k$. The class of this algebra in $\Br(k)$
is called the Clifford invariant of $\Phi$  and is denoted by $c(\Phi)$. 
In fact, $c(\Phi)$ does not change
under multiplication of $\Phi$ by an element of $k^\times$, so 
$c(\Phi)$ depends only on the projective quadric
$Q\subset \P^4_k$ with equation $\Phi=0$. Explicitly, when $\Phi$ is a diagonal 
quadratic form $\sum_{i=1}^5a_ix_i^2$, where $a_i\in k^\times$, we have 
$$c(Q)=(-1,-1)\prod_{i<j}(a_i,a_j)\in \Br(k),$$
see \cite[Chapter V, \S 3]{Lam}.

The variety parameterising projective lines contained in a smooth 3-dimensional quadric
$Q\subset\P^4_k$ is a Severi--Brauer variety 
${\rm SB}(Q)$ of dimension 3, whose class $[{\rm SB}(Q)]\in\Br(k)$ is annihilated by 2,
see \cite[Lemma 61 (1), (2)]{Kollar}. In fact, we have $[{\rm SB}(Q)]=c(Q)$, as easily follows
from  \cite[Lemma 61 (6)]{Kollar}.

Define $\SB\subset {\rm Gr}(\P^1,\P^4)\times{\Pen}$ to be the closed subscheme of pairs
$(\ell,m)$ such that the projective line $\ell\subset\P^4_k$ is contained 
in the quadric $\QQ_m=p^{-1}(m)$.
By \cite[Theorem 1.10]{R72}, the variety $\SB$ is smooth over $k$.
The projection $\pi\colon \SB\to \P^1_k$ restricts to a smooth morphism
above $\P^1_k\setminus\Sigma$ such that the fibres are Severi--Brauer varieties of dimension 3. 

As promised, we can now prove that $\delta$, as defined in terms of the singular quadrics in the pencil, lies in $(A^\times/A^{\times 2})_1$.

\ble \label{duedue}
We have 
\begin{equation}
{\rm N}_{A/k}(\delta)=\prod_{\si\in\Sigma}{\rm N}_{k_\si/k}(\delta_\si)=1
\in k^\times/k^{\times 2}, \label{Faddeev}
\end{equation}
where ${\rm N}_{A/k}\colon  A^\times/A^{\times2}\to k^\times/k^{\times 2}$ is 
induced by the norm map $A\to k$. 
\ele
{\em Proof.} We note that the morphism $\pi\colon \SB\to \P^1_k$ satisfies the assumptions
of Proposition \ref{due}. One immediately checks that
the fibre $\pi^{-1}(\si)$, $\si\in\Sigma$, consists
of two geometric components of multiplicity 1 individually defined over $k_\si(\sqrt{\delta_\si})$.
The statement thus follows from Proposition \ref{due} (i). \hfill $\Box$

\brem
{\rm
Let $X^{[2]}$ be the Hilbert scheme of $0$-dimensional subschemes of $X$ of length 2.
This is a smooth variety that can be constructed as the quotient of the blow-up of 
$X\times X$ along the diagonal ${\rm diag}(X)$,  
by the unique involution lifting the one on $X\times X$ swapping the factors. Define
$\tilde \SB\subset  X^{[2]}\times \SB$ as the set of triples $(z,\ell, m)$, where 
$z\in X^{[2]}(\bar k)$,
$(\ell,m)\in \SB(\bar k)$, and $z\subset \ell$. The projection $\tilde \SB\to \SB$
is an isomorphism outside of the 16 lines in $\SB$ given by
$(\lambda,m)$, where $\lambda\cong\P^1_{\bar k}$ is a line contained in $X$ and $m\in\Pen(\bar k)$. 
Indeed, if a line $\ell$ in $\P^4$ is contained in a quadric of $\Pen$ but not in $X$, then
intersecting $\ell$ with any quadric of $\Pen$ other than $\QQ_m$ we get a closed
subscheme $z$ of $X$ of degree 2. 
The projection $\tilde \SB\to X^{[2]}$ is an isomorphism outside of the 16 planes 
$\lambda^{[2]}=\textup{Sym}^2\lambda\cong\P^2_{\bar k}$, because then the line through $z$
is contained in a unique quadric of $\Pen$. 
 Thus both $\tilde \SB\to \SB$ and $\tilde \SB\to X^{[2]}$ 
are birational morphisms. 
}
\erem

We now study the relation between the Brauer groups of the Severi--Brauer fibration $\SB/\P^1_k$
and the Kummer fibration $Z/\P^1_k$. 
Recall that $Z$ is a smooth subvariety of $\P(\T_X)\times \Pen$.
We have a natural morphism $Z\to \P(\T_X)\to X$
whose generic fibre is a projective line. Thus $Z$ is birationally equivalent to $X\times\P^1_k$,
and the induced map $\Br(X)\to\Br(Z)$ is an isomorphism.

The natural map $Z\to Y\subset {\rm Gr}(\P^1,\P^4)\times \Pen$ gives rise to a morphism
$\varphi\colon Z\to \SB$ compatible with projections to $\Pen\cong\P^1_k$, that is,
$f=\pi\circ \varphi$.  By Theorem \ref{t1}
the fibres of $f\colon Z\to\P^1_k$ above the points of $U\subset\P^1_k$ are smooth Kummer surfaces.

\bpr \label{eins}
The homomorphism $\varphi^*\colon\Br(\SB)/\Br(k) \to \Br(Z)/\Br_0(Z)$ is zero.
\epr
{\em Proof.} 
By Proposition \ref{due} (v), to prove the statement it is enough to show that
$$\varphi^*:\H^1(k,\Pic(\SB_{\bar k}))\to \H^1(k,\Pic(Z_{\bar k}))$$ is the zero map.

Recall that if $V\to V'$ is a birational morphism between smooth projective varieties, then the induced map $\H^1(k,\Pic(V'_{\bar k}))\to \H^1(k,\Pic(V_{\bar k}))$ is an isomorphism. 
Consider the following commutative diagram, in which each variety is smooth and projective:
$$\xymatrix{Z  \ar[r]^{\varphi}\ar[d]& \SB&\tilde \SB\ar[l]\ar[d] \\
\P(\T_X) \ar[r]\ar[d]&\textup{Bl}_{{\rm diag}(X)}(X\times X)\ar[r]^{\quad \quad \ \  q}\ar[d]^p&X^{[2]}& \\
X \ar[r]^{{\rm diag} \ \ }&X\times X &&}$$
In the diagram, ${\rm diag}$ is the diagonal embedding, $\textup{Bl}_{{\rm diag}(X)}(X\times X)$ is the blow-up of $X\times X$ along the diagonal,  the map $\P(\T_X) \to \textup{Bl}_{{\rm diag}(X)}(X\times X)$ is the inclusion of the exceptional divisor of the blow-up, and the maps $p$ and $q$ are the natural 
contraction and quotient maps respectively. 
Each vertical arrow in the diagram is a birational morphism, with the exception of the projection $\P(\T_X) \to X$ which nevertheless induces an isomorphism between 
$\H^1(k,\Pic(X_{\bar k}))$ and $\H^1(k,\Pic(\P(\T_X)_{\bar{k}}))$.

From consideration of the diagram, it suffices to show that the map 
\[{\rm diag}^*\circ (p^*)^{-1}\circ q^*:\H^1(k,\Pic(X^{[2]}_{\bar{k}})) \to \H^1(k,\Pic(X_{\bar{k}}))\]
is trivial. Since $X_{\bar{k}}$ is rational, the two projections $X\times X \to X$ induce an isomorphism $\Pic(X_{\bar k}\times X_{\bar{k}})\cong \Pic(X_{\bar{k}})\oplus \Pic(X_{\bar{k}})$, after which 
the map ${\rm diag}^*$ identifies with the sum map $\Pic(X_{\bar{k}})\oplus \Pic(X_{\bar{k}})\to \Pic(X_{\bar{k}})$. Since $2\H^1(k,\Pic(X_{\bar k}))=0$
(Proposition \ref{coh dp}), it suffices to show that the image of the map 
\[(p^*)^{-1}\circ q^*: \H^1(k,\Pic(X^{[2]}_{\bar{k}}))\to  
\H^1(k,\Pic(X_{\bar k}\times X_{\bar{k}}))\] 
is invariant under the action of the involution $\iota$ of $X\times X$ swapping the factors. This is a formal consequence of the fact that $X^{[2]}$ is the quotient of $\textup{Bl}_{{\rm diag}(X)}(X\times X)$ by an involution lifting $\iota$. \hfill $\Box$

\medskip

Let $\Br_{\rm vert}(Z)$ be the vertical
Brauer group associated to the morphism $f\colon Z\to \P^1_k$.

\bco
The natural map $\Br(k)\to\Br_{\rm vert}(Z)$ is surjective.
\eco
{\em Proof.} The vertical Brauer group $\Br_{\rm vert}(Z)$ is the subgroup of 
$\Br(Z_{k(t)})$ defined
as the intersection of $\Br(Z)$ with the image of $\Br(k(t))\to \Br(Z_{k(t)})$. 
The morphism $\varphi\colon Z\to \SB$ is compatible with the projections $f$ and $\pi$ 
to $\P^1_k$,
hence it induces a homomorphism $\varphi^*\colon \Br_{\rm vert}(\SB)\to \Br_{\rm vert}(Z)$.
Just like the fibres of $\pi\colon \SB\to\P^1_k$, the fibres of $f$ above 
$\P^1_k\setminus\Sigma$ are split 
and the fibre $Z_\si$, $\si\in\Sigma$, decomposes as the union of 
two geometrically irreducible components individually defined
over $k_\si(\sqrt{\delta_\si})$. This implies that the map
$\varphi^*\colon \Br_{\rm vert}(\SB)\to \Br_{\rm vert}(Z)$ is surjective.
Since the composition
$$\Br_{\rm vert}(\SB)/\Br(k)\to \Br_{\rm vert}(Z)/\Br_0(Z)\hookrightarrow\Br(Z)/\Br_0(Z)$$
is trivial by Proposition \ref{eins}, the result follows. \hfill $\Box$

\brem  \label{Brauer of X rem}
{\rm In spite of Proposition \ref{eins}, there are canonical isomorphisms
$$\H^1(k,\Pic(Z_{\bar k}))\cong\H^1(k,\Pic(X_{\bar k}))\cong
 {\rm B}_{A,\delta}\cong \H^1(k,\Pic(\SB_{\bar k})) \cong \Br(\SB)/\Br(k).$$
The first one is a consequence of stably birational equivalence of $X$ and $Z$, and the rest 
follow from the combination of Propositions \ref{due} and \ref{coh dp}.}
\erem

\section{Main results} \label{S3}

\subsection{Admissible local conditions}

From now on we let $k$ be a number field. 

Our method proceeds by choosing a parameter $b\in k$, $P(t)\neq 0$, satisfying certain
local conditions made precise in the following definition.
Recall that we have a finite \'etale degree 2 morphism of $0$-dimensional $k$-schemes
$\Delta\to\Sigma$.
The Galois group $\Ga$ acts on $\Delta(\bar k)$ via its quotient $\Gal(k_\delta/k)$.
The $k$-scheme $\Delta$ is a disjoint union of closed points which bijectively correspond
to the orbits of $\Gal(k_\delta/k)$ in $\Delta(\bar k)=\Delta(k_\delta)$.

\bde {\rm 
A pair $(C,D)$, where $C$ is a conjugacy class in 
$\Gal(k_\delta/k)$ and $D$ is a closed point of $\Delta$
such that any $g\in C$ has a fixed point in $D(\bar k)$, will be called an
{\em admissible local condition} for $(P(t),\delta)$ over $k$.
We say that $b\in k$, $P(b)\neq 0$, satisfies an admissible local condition
$(C,D)$ at a place $w$ of $k$, if we have
\begin{itemize}
\item $w\notin S_0$;
\item ${\rm Frob}_w\in C$;
\item $\val_w(Q(b))=1$, where $Q(t)$ is the monic irreducible factor of $P(t)$ such that the image of 
$D$ in $\Sigma$ is given by $Q(t)=0$;
\item $D$ has a $k_w$-point over the $k_w$-point of $\Sigma$ given by 
the root $\theta_w\in O_w$ of $Q(t)$ such that $\theta_w-b$ is in the maximal ideal of $O_w$.
\end{itemize}
}
\ede

\brem\label{realise}{\rm
For any admissible local conditions $(C_i,D_i)$, $i=1,\ldots, n$, for $(P(t),\delta)$ over $k$
there exist $b\in k$, $P(b)\neq 0$, and places $v_1,\ldots,v_n$ of $k$
such that $b$ satisfies $(C_i,D_i)$ at $v_i$ for $i=1,\ldots, n$.
Indeed, by the Chebotarev density theorem we can find $v_1,\ldots,v_n$ not in $S_0$
such that  for $i=1,\ldots,n$ we have ${\rm Frob}_{v_i}\in C_i$.
Let $P_i(t)$ be the monic irreducible factor of $P(t)$ such that the image of 
$D_i$ under the map $\Delta\to\Sigma$ is given by $P_i(t)=0$.
Admissibility implies that for $i=1,\ldots,n$ the polynomial
$P_i(t)$ has a simple root $\theta_i\in O_{v_i}$,
moreover, $\Delta$ has a $k_{v_i}$-point over $\theta_w\in\Sigma(k_w)$. Let
$b_i\in O_{v_i}$ be such that $\val_{v_i}(b_i-\theta_i)=1$. 
Then $P_i(b_i)\neq 0$ and $\val_{v_i}(P_i(b_i))=1$.
By weak approximation, we can find an element $b\in k$ close enough to $b_i$ in 
$k_{v_i}$ for $i=1,\ldots,n$. Then $b$ satisfies $(C_i,D_i)$ at $v_i$ for $i=1,\ldots, n$.
}\erem

Let $C_b$ be the curve of genus 2 with affine equation 
$y^2=(x-b)P(x)$, where $b\in k$, $P(b)\neq 0$, and let $J_b$ be the Jacobian of $C_b$. 
If $F$ is a quadratic extension of $k$, or $F=k$, we denote by $J_b^F$ the quadratic twist of 
$J_b$ by $F$. Clearly, $J_b^F$ is the Jacobian of the curve $y^2=a(x-b)P(x)$, where
$F=k(\sqrt{a})$.

\ble \label{mumbai3}
Let $w$ be an odd prime. Suppose that the group $\F_w$-scheme of connected components
of the special fibre of the N\'eron model of $J_b^F$ over $O_w$ is a constant abelian group $\pi_0$.
Then there is an exact sequence
\begin{equation}
0\to \partial(J_b^F(k_w))\cap \H^1_\nr(k_w,G)\to \H^1_\nr(k_w,G)\to\pi_0[2]\to 0.
\label{mumbai2}
\end{equation}
\ele
{\em Proof.}
We have a commutative diagram with exact rows
$$\xymatrix{
0\ar[r]&J_b^F(k_w^\nr)/2\ar[r]^\partial&\H^1(k_w^\nr,G)\ar[r]&\H^1(k_w^\nr, J_b^F)[2]\ar[r]&0\\
0\ar[r]&J_b^F(k_w)/2\ar[r]^\partial\ar[u]&\H^1(k_w,G)\ar[r]\ar[u]&\H^1(k_w, J_b^F)[2]\ar[r]\ar[u]&0}$$
where $\partial$ is the connecting map defined by the exact sequence of $\Ga$-modules
$$0\to G\to J_b^F\xrightarrow{[2]}J_b^F\to 0.$$ The snake lemma gives rise to the exact sequence
\begin{equation}
0\to \partial(J_b^F(k_w))\cap \H^1_\nr(k_w,G)\to \H^1_\nr(k_w,G)\to \H^1_\nr(k_w,J^F_b)[2],
\label{mumbai1}
\end{equation}
where $\H^1_\nr(k_w,J^F_b)$ is by definition the kernel of 
$\H^1(k_w, J_b^F)\to \H^1(k_w^\nr, J_b^F)$.
By \cite[Proposition I.3.8]{Mil06}, we have 
$\H^1_\nr(k_w,J^F_b)\cong\H^1(\F_w,\pi_0)$.
Thus the last arrow in (\ref{mumbai1}) is the map 
$J^F_b(k_w^\nr)[2]/(\Frob_w-1)\to \pi_0[2]$.

The natural reduction map $J^F_b(k_w^{\textup{nr}})\to \pi_0$ is surjective by Hensel's lemma.
The kernel of this map is ${\mathcal J}^0(O_w^\nr)$, where ${\mathcal J}^0$ is
the identity component  of the N\'eron model ${\mathcal J}$ of $J^F_b$ over $O_w$.
The morphism $[2]\colon{\mathcal J}^0\to{\mathcal J}^0$ is 
\'etale and surjective since $w$ is an odd prime, thus 
${\mathcal J}^0(O_w^\nr)$ is a $2$-divisible  abelian group.
Applying the snake lemma to the exact sequence
$$0\to {\mathcal J}^0(O_w^\nr)\to J^F_b(k_w^{\textup{nr}})\to \pi_0\to 0$$
we obtain the surjectivity of the reduction map $J^F_b(k_w^\nr)[2]\to \pi_0[2]$.
\hfill $\Box$

\medskip

 Write $J^F_{b,\delta}$ for the $k$-torsor obtained by twisting $J^F_b$
by the 1-cocycle $\tilde\delta\colon\Ga\to G(\bar k)$ with respect to the action of 
$G\cong J^F_b[2]$ by translations.
Recall that $\Kum(J_{b,\delta})$ is the Kummer surface attached to the 2-covering $J_{b,\delta}$.
We have $\Kum(J_{b,\delta})\cong \Kum(J^F_{b,\delta})$ for any $F$.

The next proposition explains the relevance of admissible local conditions.

\bpr \label{jan1}
Let $b\in k$, $P(t)\neq 0$, and let $w$ be a place of $k$ outside of $S_0$ 
such that $\val_w(P(b))=1$. The following conditions are equivalent.

{\rm (i)} $b$ satisfies an admissible local condition at $w$.

{\rm (ii)} $J^{F_w}_{b,\delta}(k_w)\neq\emptyset$ where $F_w$ is an unramified extension $k_w$
of degree at most $2$.

{\rm (iii)} $\Kum(J_{b,\delta})(k_w)\neq\emptyset$.
\epr
{\em Proof.} The conditions on $w$ imply that the equation $y^2=(x-b)P(x),$
 viewed over  $O_w$ and considered along with the usual chart at infinity,  defines a semistable model of $C_b$ whose special fibre has a single node at $x=b$. The minimal proper regular model of $C_b$ over $O_w$ is obtained by blowing up once at the node. Let $\mathcal{J}_b\to \textup{Spec}(O_w)$ be the N\'{e}ron model of $J_b$. Since the special fibre of the minimal regular model 
of $C_b$ consists of two components intersecting in two ordinary double points, 
\cite[Proposition 9.6.10]{BLR90} shows that 
the group $\F_w$-scheme of connected components of 
$\mathcal{J}_b\times_{O_w}\F_w$ is $\pi_0\cong\Z/2$. 
By Lemma \ref{mumbai3} we have an exact sequence
\begin{equation}
0\to \partial(J_b(k_w))\cap \H^1_\nr(k_w,G)\to \H^1_\nr(k_w,G)\to \Z/2\to 0,
\label{mumbai4}
\end{equation}
where the surjective map is $G(k^\nr_w)/(\Frob_w-1)\to \pi_0/(\Frob_w-1)$
induced by a natural map $s\colon G(k^\nr_w)\to \pi_0$. Let us compute $s$.
By assumption, there is a root $\theta_w\in O_w$ of $P(x)$
such that $\val_w(b-\theta_w)=1$. 
The difference of Weierstrass points $(b,0)-(\theta_w,0)$ on $C_b$
is a point in $J_b(k_w)[2]$ whose reduction belongs to the non-identity component of 
$\mathcal{J}_b\times_{O_w}\F_w$, so $s$ sends this point to
$1\in\Z/2$. If $\xi\in k^\nr_w$ is a root of $P(t)=0$,
$\xi\neq \theta_w$, then the reduction of $(b,0)-(\xi,0)\in J_b(k^\nr_w)[2]$ belongs to the identity
component of $\mathcal{J}_b\times_{O_w}\F_w$, so $s$ sends this point to $0\in\Z/2$.
Thus the map $s$ is the composition of the natural embedding $G\to R_{A/k}(\mu_2)$
with the projection $R_{A/k}(\mu_2)\to\mu_2$ given by the factor $x-\theta_w\in k_w[x]$ of $P(t)$.
This implies that the surjective map in (\ref{mumbai4}) sends $\delta\in \H^1_\nr(k_w,G)$ 
to 0 or 1 in $\Z/2$ if the preimage in $\Delta$ of the $k_w$-point of $\Sigma$ given by $x=\theta_w$ 
consists of two $k_w$-points or one closed point, respectively.

Let us prove that (i) and (ii) are equivalent. We note that $b$ satisfies an admissible
local condition at $w$ if and only if the preimage in $\Delta$ of the $k_w$-point of 
$\Sigma$ given by $x=\theta_w$ consists of two $k_w$-points. Thus (i) is equivalent to
$\delta\in \partial(J_b^F(k_w))$, which is also equivalent to $J_{b,\delta}(k_w)\neq\emptyset$.

N\'{e}ron models commute with unramified base change, so the same considerations
apply to $J_b^{F_w}$. 

Let us prove that (i) and (iii) are equivalent.
We have $Z_b=\Kum(J_{b,\delta})$. Since $\val_w(P(b))=1$, the polynomial
$P(x)$ has a simple root $\theta_w\in O_w$. The surface $Z_b$ over $k$
naturally extends to a projective scheme ${\mathcal Z}_b\to \Spec(O_w)$
with special fibre isomorphic to the reduction of the singular fibre $Z_{\theta_w}$
of $f\colon Z\to\P^1_k$.
Moreover, since smoothness is stable under base change, 
the morphism ${\mathcal Z}_b\to \Spec(O_w)$ is smooth away from the singular
locus of the special fibre. As discussed at the end of Section \ref{one},
the scheme over $\Sigma$ parameterising the geometrically irreducible components
of the singular fibres of $f\colon Z\to\P^1_k$ is isomorphic to $\Delta$.
Thus the reduction of $Z_{\theta_w}$ is the union of
two irreducible components if and only if ${\rm Frob}_w$ fixes a point of $\Delta(k^\nr_w)$
above the root $\theta_w$ of $P(x)$.

If this is the case, each component is the reduction of $X$ blown up at four $\bar k$-points,
meeting the other component in a smooth curve of genus 1.
By the definition of $S_0$, there are $\F_w$-points in $X$ outside this curve. Thus
the reduction of $Z_{\theta_w}$ has a smooth $\F_w$-point.
By Hensel's lemma it can be lifted by a point in $Z_b(k_w)$. 

In the opposite case,
the reduction of $Z_{\theta_w}$ is not split so that $Z_b(k_w)=\emptyset$ by
\cite[Corollary 10.1.9]{CTS21} together with the valuative criterion of properness. \hfill $\Box$

\medskip

The equivalence of (i) and (ii) can be proved in the same way as the equivalence of (i) and (iii).
Consider the N\'eron model over $\P^1_k$ of the Jacobian $J_t$
of the genus 2 curve over $k(t)$ given by $y^2=(x-t)P(x)$.
Arguing as in the beginning of the proof we see that the singular fibres above the closed points
of $\Sigma$ consist of two irreducible components. Twisting by the 1-cocycle $\tilde\delta\colon
\Gamma\to G$ (which is unramified since $w\notin S_0$)
we obtain a model of the torsor $J_{t,\delta}$. The scheme over $\Sigma$
parameterising the geometrically irreducible components of this fibration is $\Delta$, just like
for the Kummer fibration $f\colon Z\to\P^1_k$. The rest of the proof is similar to 
the proof of the equivalence of (i) and (iii), using the N\'eron property of rational points.

\subsection{Hasse principle: from Kummer to del Pezzo}

Let $\K_{b,\delta}$ be the fibre of the Kummer fibration $\K_\delta\to U$ at $b\in k$, $P(b)\neq 0$.
Thus $\K_{b,\delta}=\Kum(J_{b,\delta})$.
Using notation of \cite{SZ17}, we denote by $\Pi_\delta$ the Kummer lattice in 
$\Pic(\K_{b,\delta,\bar k})$, i.e.~the saturation
of the free abelian group freely generated by the classes of the 16 lines.
The group $G$ naturally acts on $\Pi_\delta$, as does the Galois group $\Ga$. 
Neither action depends on $b$.  

\bthe \label{dp}
Let $X$ be a quartic del Pezzo surface over a number field $k$ given by a pair $(P(t),\delta)$.
Then $X(\bA_k)^\Br\neq\emptyset$ implies $X(k)\neq\emptyset$, provided
the following holds:

There exist finitely many admissible local conditions for $(P(t),\delta)$ over $k$
such that, for any finite extension $K/k$ linearly disjoint from $k_\delta/k$
and any $b\in K$, $P(b)\neq 0$, satisfying these local conditions over $K$,
the Kummer surfaces $\K_{b,\delta}=\Kum(J_{ b,\delta})$
have the property that the algebraic Brauer--Manin obstruction
attached to the Kummer lattice $\Pi_\delta\subset\Pic(\Kumm_{b, \delta, \bar k})$
is the only obstruction to the Hasse principle on $\Kumm_{b,\delta}$.   
\ethe
{\em Proof.} Let $X$ be a quartic del Pezzo surface given by $(P(t),\delta)$. 
Then $X$ is stably birational to $Z$, and the restriction of $f\colon Z\to\P^1_k$
to $U$ is isomorphic to the Kummer fibration $\K_\delta\to U$, 
see Theorem \ref{t1}. In particular,
the morphism $Z_U\to U$ is smooth and projective with geometrically integral fibres.
This implies that the restriction to the generic fibre induces an isomorphism
$\Pic(Z_{U_{\bar k}})\tilde\lra\Pic(Z_{\bar k(t)})$.
In particular, $\Pic(Z_{U_{\bar k}})$ is a finitely generated free abelian group,
hence $\H^1(k,\Pic(Z_{U_{\bar k}}))$ is finite.  

The proof of the following lemma is well-known and is given here for the convenience of the reader.

\ble \label{lem:Br_to_pic_surj}
The natural map 
$\Br_1(Z_U)\to \H^1(k,\Pic(Z_{U_{\bar k}})) $
is surjective. 
\ele
{\em Proof.} The Leray spectral sequence 
$\H^p(k,\H^q(Z_{U_{\bar k}},\G_m))\Rightarrow \H^{p+q}(Z_U,\G_m)$
gives an exact sequence
$$\Br_1(Z_U)\to \H^1(k,\Pic(Z_{U_{\bar k}})) \to \H^3(k,\H^0(Z_{U_{\bar k}},\G_m)),$$
so it suffices to show that $\H^3(k,\H^0(Z_{U_{\bar k}},\G_m))=0$. Since 
the morphism $Z_{U_{\bar k}}\to U_{\bar k}$ is proper with geometrically integral fibres,
we have $\H^0(Z_{U_{\bar k}},\G_m)\cong \H^0(U_{\bar k},\G_m)\cong \bar k[U]^\times$.
There is an exact sequence of $\Ga$-modules
$$0\to \bar k^\times\to \bar k[U]^\times \to \Z[\Sigma]\to 0,$$
where the third arrow sends an invertible regular function on 
$U_{\bar k}=\A^1_{\bar k}\setminus \Sigma_{\bar k}$ to its divisor in 
$\A^1_{\bar k}$. We obtain an exact sequence of Galois cohomology groups
$$\H^3(k,\bar k^\times)\to \H^3(k,\bar k[U]^\times) \to \H^3(k,\Z[\Sigma]).$$ 
Since $k$ is a number field, we have $\H^3(k,\bar k^\times)=0$ \cite[Ch.~VII, \S 11.4]{CF}.
We also have $\H^3(k,\Z[\Sigma])=0$, which is a consequence of Shapiro's lemma and the 
vanishing of  $\H^3(k',\Z)\cong\H^2(k',\Q/\Z)$ for any number field $k'$   
\cite[\S 4, Satz 1 (ii)]{Sch79}. \hfill $\Box$

\medskip

The generic fibre $Z_{\bar k(t)}$ is the Kummer surface attached to the Jacobian 
of a genus 2 curve all of whose Weierstrass points are $\bar k(t)$-rational. 
Thus the 16 lines on $Z_{\bar k(t)}$ are defined over $\bar k(t)$. Since $Z_{\bar k(t)}$
has a $\bar k(t)$-point, 
we have $\Pic(Z_{\bar k(t)})=\Pic(Z_{\ov{k(t)}})^{\Gal(\ov{k(t)}/\bar k(t))}$,
see \cite[Remark 5.4.3 (1)]{CTS21}. In particular, $\Pic(Z_{\bar k(t)})$ is a saturated subgroup
of $\Pic(Z_{\ov{k(t)}})$. This implies that the Kummer lattice of $\Pic(Z_{\ov{k(t)}})$
is contained in $\Pic(Z_{\bar k(t)})\cong \Pic(Z_{U_{\bar k}})$.
Thus we have natural embeddings of $\Ga$-modules $\Pi_\delta\to\Pic(Z_{U_{\bar k}})$
and $\Pi_\delta\to\Pic(Z_{b,\bar k})$ compatible under the specialisation map
$\Pic(Z_{U_{\bar k}})\to\Pic(Z_{b, \bar k})$. 

By Lemma \ref{lem:Br_to_pic_surj}, we can choose lifts to $\Br_1(Z_U)$ of each element 
of the image of $\H^1(k,\Pi_\delta)$ in $\H^1(k,\Pic(Z_{U_{\bar k}}))$. 
Define $B$ as the subgroup of $\Br_1(Z_U)$ generated by these lifts.
Thus $B$ is a finite subgroup of $\Br_1(Z_U)$ such that for any $b\in k$, $P(b)\neq 0$,
the image of the specialisation map 
$$B\hookrightarrow\Br_1(Z_U)\to \Br_1(Z_b)\to \H^1(k,\Pic(Z_{b, \bar k}))$$
contains the image of $\H^1(k,\Pi_\delta)\to\H^1(k,\Pic(Z_{b, \bar k}))$.

The variety $Z$ is stably birational to $X$,
hence if $X$ is everywhere locally soluble, then so is $Z$.
For the same reason, the morphism $Z\to X$
induces an isomorphism $\Br(X)\tilde\lra\Br(Z)$, thus if $X(\bA_k)^\Br\neq\emptyset$, then $Z(\bA_k)^\Br\neq\emptyset$, see \cite[Proposition 13.3.11]{CTS21}.
 
For any sufficiently large finite set of places $S$, the projective $\P^1_k$-scheme $Z$ extends to  a projective $\P^1_{O_{k,S}}$-scheme ${\mathcal Z}$ such that the morphism ${\mathcal Z}\to \P^1_{O_{k,S}}$ is smooth away from $P(t)=0$.
We choose $S$ such that $S_0\subset S$ and the Brauer classes 
$\alpha_1,\ldots,\alpha_r$ from the proof of Proposition \ref{h-w} have zero evaluation maps outside of $S$.

Using Remark \ref{realise} we find $b\in k$, $P(b)\neq 0$, 
satisfying the admissible local conditions 
in the statement of the theorem over $k$. Then 
$Z_b(k_{v_i})\neq\emptyset$ for $i=1,\ldots,n$ by Proposition \ref{jan1}.

An application of Proposition \ref{h-w} gives a closed point $R=\Spec(K)$ in $U$ such that the fibre
$Z_R$ is a Kummer surface over $K$ with an adelic point orthogonal to $B$.
We can ensure that $[K:k]$ is coprime to $30$, hence $K$ and $k_\delta$
are linearly disjoint over $k$. Moreover, we can ensure that $K$ has a place $w_i$ 
over $v_i$ such that $K_{w_i}\cong k_{v_i}$ and $R$ is
arbitrarily close to $b$ at $w_i$, for $i=1,\ldots, n$. It follows that 
$R\in \A^1_k(K)$ satisfies the local conditions in the statement of the theorem over $K$.

The condition in the theorem now implies that $Z_R$ has a $K$-point. It follows that $X$
has a $K$-point. Since $[K:k]$ is odd,
the Amer--Brumer theorem \cite{A76, B78} then implies that $X$ has a $k$-point.
\hfill $\Box$

\subsection{Proof of Theorem A}

We shall use the following restatement of the main result of \cite{Har19}.

\bpr[Harpaz] \label{h-thm}
Let $k$ be a number field. Assume that $P(t)$ is completely split over $k$, so that
$P(t)=\prod_{i=1}^5(t-a_i),$ where $a_i\in k$, $a_i\neq a_j$ for $i\neq j$.
Assume that $\delta=(\delta_1,\ldots,\delta_5)\in (k^\times/k^{\times 2})^5$ is such that
$\prod_{i=1}^5\delta_i=1$ is the only multiplicative relation satisfied by the $\delta_i$. 
Let $b\in k$ be such that $P(b)\neq 0$. 
Assume that $\Sha(J^F_b)[2^\infty]$ is finite for every extension $F/k$ of degree at most $2$.
Then the Kummer surfaces $\K_{b,\delta}=\Kum(J_{b,\delta})$
have the property that the algebraic Brauer--Manin obstruction
attached to the Kummer lattice $\Pi_\delta$
is the only obstruction to the Hasse principle, provided 
$b$ satisfies admissible local conditions at places
$w_1,\ldots,w_5$ of $k$ outside of $S_0$ such that, for $i=1,\ldots,5$, we have
$\val_{w_i}(b-a_i)=1$ and the image of $\delta_i$ in 
$k_{w_i}^\times/k_{w_i}^{\times 2}$ is trivial, whereas the image of $\delta_j$ 
for some $j\neq i$ is not.
\epr
{\em Proof.} From \cite[Theorem 1.3, Remark 2.9]{Har19} we know that under the 
assumptions of the proposition
the algebraic Brauer--Manin obstruction is the only obstruction to the Hasse principle on
$\Kum(J_{b,\delta})$. In fact, Harpaz proves the stronger statement concerning
the algebraic Brauer--Manin obstruction attached to the Kummer lattice $\Pi_\delta$.
Let us explain this. Write $\Kumm_{b,\delta}=\Kum(J_{b,\delta})$ and let 
$W=\Kumm_{b,\delta}\setminus \Lambda_\delta$
be the complement to the 16 lines. 
 The diagram (7) from \cite{SZ17} is
a commutative diagram of $\Ga$-modules with exact rows and columns
\begin{equation}\begin{array}{ccccccccc}
&&&&0&&0&&\\
&&&&\uparrow&&\uparrow&&\\
&&&&\NS(J_{b,\delta, \bar k})&=&\NS(J_{b,\delta, \bar k})&&\\
&&&&\uparrow&&\uparrow&&\\
0&\lra& \Z[\Lambda_\delta]&\lra& \Pic(\Kumm_{b,\delta,\bar{k}}) &\lra &\Pic(W_{\bar k})&\lra& 0\\
&&||&&\uparrow&&\uparrow&&\\
0&\lra& \Z[\Lambda_\delta]&\lra& \Pi_\delta &\lra& \Pic(W_{\bar k})_{\rm tors}&\lra& 0\\
&&&&\uparrow&&\uparrow&&\\
&&&&0&&0&&
\end{array}\label{u4}
\end{equation}
The $\Ga$-module $\Z[\Lambda_\delta]$ is permutational, hence $\H^1(k,\Z[\Lambda_\delta])=0$.
Thus (\ref{u4}) gives rise to a commutative diagram with exact columns
\begin{equation}\begin{array}{ccccccccc}
\H^1(k,\NS(J_{b,\delta, \bar k}))&=&\H^1(k,\NS(J_{b,\delta, \bar k}))\\
\uparrow&&\uparrow\\
\H^1(k, \Pic(\Kumm_{b,\delta,\bar{k}})) &\hookrightarrow &\H^1(k,\Pic(W_{\bar k}))\\
\uparrow&&\uparrow\\
\H^1(k, \Pi_\delta) &\hookrightarrow& \H^1(k,\Pic(W_{\bar k})_{\rm tors})&\\
\end{array}\label{u5}
\end{equation}
From this diagram we conclude that the subgroup of $\H^1(k, \Pic(\Kumm_{b,\delta,\bar{k}}))$
consisting of the elements whose images in $\H^1(k,\Pic(W_{\bar k}))$ come from
$\H^1(k,\Pic(W_{\bar k})_{\rm tors})$ is just the image of $\H^1(k, \Pi_\delta)$. 
It follows that the group denoted by 
$${\mathcal C}(\Kumm_{b,\delta})\subset\H^1(k, \Pic(\Kumm_{b,\delta,\bar{k}}))=\Br_1(\Kumm_{b,\delta})/\Br_0(\Kumm_{b,\delta})$$ 
in \cite[Definition 4.3]{Har19} is the image of $\H^1(k,\Pi_\delta)$. 
By \cite[Remark 2.9]{Har19} the obstruction given by ${\mathcal C}(\Kumm_{b,\delta})$
is the only obstruction to the Hasse principle on $\Kumm_{b,\delta}$.~\hfill $\Box$

\medskip

The following theorem gives the Hasse principle for Kummer surfaces
with irreducible polynomial $P(t)$, assuming that $b\in k$ satisfies 
certain local conditions. 
The proof will be given in Sections \ref{SAF} and \ref{D10}.

\bthe \label{thm}
Let $k$ be a number field. Assume that a pair $(P(t),\delta)$ is such that $P(t)$ is irreducible.
Let $b\in k$ be such that $P(b)\neq 0$. 
Assume that $\Sha(J^F_b)[2^\infty]$ is finite for every extension $F/k$ of degree at most $2$.

{\rm (a)} Suppose that the Galois group of $P(t)$ is ${\rm S}_5$, $\AA_5$, or $\FF_{20}$. 
Then the Hasse principle holds 
for the Kummer surfaces $\Kum(J_{b,\delta})$ for which there is a
prime $w$ outside of $S_0$ such that $\val_w(P(b))=1$.

{\rm (b)} Suppose that the Galois group of $P(t)$ is $\DD_{10}$.
Then the Hasse principle holds for the Kummer surfaces $\Kum(J_{b,\delta})$
for which there exist primes $w_1$ and $w_2$ outside of $S_0$ such that 
for $i=1,2$ we have $\val_{w_i}(P(b))=1$,
the Frobenius element at $w_i$ acts on the roots of $P(t)$ as a double transposition, and
$\loc_{w_1}(\delta)\neq 0$, $\loc_{w_2}(\delta)=0$. 
\ethe

Now we are ready to prove Theorem A of the Introduction.

\bthe \label{adrasan2}
Let $k$ be a number field. Assume that $\Sha({\rm Jac}(C))[2^\infty]$ is finite
for every genus $2$ curve $C$ over every finite extension of $k$. 
Then the quartic del Pezzo surfaces $X$ over $k$ such that the polynomial 
$P(t)\in k[t]$ is completely split and $\Br(X)=\Br_0(X)$, or $P(t)\in k[t]$ is
irreducible, satisfy the Hasse principle.
\ethe
{\em Proof.} Without loss of generality we can assume that $\delta\neq 0$.
Indeed, in the opposite case $X$ contains a $k$-line, so there is nothing to prove.

Suppose that $P(t)\in k[t]$ is completely split, so that
$P(t)=\prod_{i=1}^5(t-a_i)$, where $a_i\in k$, $a_i\neq a_j$ when $i\neq j$.
Suppose also that $\Br(X)=\Br_0(X)$. This implies that
$\Gal(k_\delta/k)$ is a subgroup of $(\Z/2)^4\rtimes {\rm S}_5$ contained in $(\Z/2)^4$.
We can assume without loss of generality
that $\Gal(k_\delta/k)$ does not fix a point of $\Delta(\bar k)$ or, equivalently, that
no $\delta_i$ is $1\in k^\times/k^{\times 2}$.
Indeed, otherwise $X$ has a conic bundle defined over $k$, so
there is a morphism $X\to \P^1_k$ making $X$ a conic bundle with 4 geometric degenerate fibres. 
In this case $X$ satisfies the Hasse principle unconditionally by a theorem of Salberger \cite[Corollary 7.8]{Sal88}. 

We claim that $\Gal(k_\delta/k)$ is the normal subgroup 
$(\Z/2)^4$ of $(\Z/2)^4\rtimes {\rm S}_5$.
Indeed, if $\Gal(k_\delta/k)$ is a proper subgroup of $(\Z/2)^4$ and no 
$\delta_i$ is $1\in k^\times/k^{\times 2}$, then ${\rm B}_{A,\delta}\neq 0$. 
By Proposition \ref{coh dp}
we have $\H^1(k,\Pic(X_{\bar k}))\cong{\rm B}_{A,\delta}\neq 0$. Since $k$ is a number field,
we have $\H^3(k,\bar k^\times)=0$. Then the exact sequence of low degree terms of the 
spectral sequence $\H^p(k,\H^q(X_{\bar k},\G_m))\Rightarrow\H^{p+q}(X,\G_m)$
gives an isomorphism
$\Br(X)/\Br_0(X)\cong \H^1(k,\Pic(X_{\bar k}))$, hence we get a contradiction.

For $i=1,\ldots,5$ take $g_i$ to be the product of all generators $c_j\in(\Z/2)^5$ except $c_i$,
and let $D_i$ be the closed point of $\Delta$ above the $k$-point of $\Sigma$ given by $t=a_i$.
The local conditions $(g_i, D_i)$ are admissible.
Let $w_1,\ldots,w_5$ be primes outside of $S_0$ such that ${\rm Frob}_{w_i}=g_i\in(\Z/2)^4$.
By Lemma \ref{delta zero} we have $\loc_{w_i}(\delta)\neq 0$.
Then Proposition \ref{h-thm}
implies that the assumption of Theorem \ref{dp} is satisfied.

Assume now that $P(t)$ is irreducible over $k$.

Let $\Gal(P)$ be one of ${\rm S}_5$, $\AA_5$, $\FF_{20}$.
By Theorem \ref{thm} (a), the Hasse principle holds 
for the Kummer surfaces $\Kum(J_{b,\delta})$ for which there is a
prime $w$ outside of $S_0$ such that $\val_w(P(b))=1$.
In particular, this holds if $b$ satisfies an admissible local condition at $w$. 
Thus the assumption of Theorem \ref{dp} is satisfied with one arbitrary admissible local condition.

Let $\Gal(P)\cong\DD_{10}$. 
By Theorem \ref{thm} (b) the assumption of Theorem \ref{dp} is satisfied
for the admissible local conditions at two primes $w_1, w_2\notin S_0$ whose Frobenius elements
in $(\Z/2)^4\rtimes \DD_{10}$ are $g_1=c_1c_3(12)(34)$ and $g_2=(12)(34)$, respectively.
Indeed, Lemma \ref{delta zero} implies that $\loc_{w_1}(\delta)\neq 0$, 
$\loc_{w_2}(\delta)=0$, as required.

Finally, let $\Gal(P)\cong\Z/5$. From Lemma \ref{is_iso_vart} 
the non-triviality of $\delta$ implies that $\Gal(k_\delta/k)$ coincides with
$(\Z/2)^4\rtimes\Z/5$ in $(\Z/2)^4\rtimes {\rm S}_5$. Let $K$ be the splitting field of $P(t)$.
Then $X_K$ is a quartic del Pezzo surface such that $P(t)$ is completely split over $K$
and $\Gal(K_\delta/K)$ is the subgroup $(\Z/2)^4$ of $(\Z/2)^4\rtimes {\rm S}_5$. 
This implies $\Br(X_K)=\Br_0(X_K)$. Applying an already proved case of the theorem, we obtain
that $X$ has a $K$-point. Since $[K:k]=5$ is odd, 
the Amer--Brumer theorem \cite{A76, B78} implies that $X$ has a $k$-point. \hfill $\Box$

\medskip

\brem{\rm The argument in the last paragraph of the above proof works equally well for
quartic del Pezzo surfaces such that $\Gal(k_\delta/k)=(\Z/2)^4\rtimes\Z/3$.}
\erem

\brem \label{split}{\rm
For quartic del Pezzo surfaces $X$ over $k=\Q$ such that
the polynomial $P(t)\in k[t]$ is completely split and $\Br(X)=\Br_0(X)$, 
one can prove Theorem \ref{adrasan2} using the fibration theorem of Harpaz--Wittenberg
for rational points \cite[Theorem 9.17]{HW16} (see also \cite[Theorem 14.2.21]{CTS21}). 
Indeed, since $P(t)$ completely splits, the complement
$\P^1_k\setminus U$ is a union of $k$-points, so that Hypothesis (HW) holds for $k=\Q$,
see \cite[Theorem 9.14]{HW16} based on work of L.~Matthiesen 
in additive combinatorics \cite{Mat18}. Thus \cite[Theorem 9.17]{HW16}
implies the existence of a point $b\in k$, $P(b)\neq 0$, such that $Z_{b}$ has an adelic point
orthogonal to $B$. 
This proof avoids replacing $k$ by a finite extension. 
As a consequence, the assumption that $\Sha(J^F_b)[2^\infty]$ 
is finite for every $b\in k$, $P(b)\neq 0$, and
every extension $F/k$ of degree at most $2$, is enough to conclude that
$X$ satisfies the Hasse principle. 
When $\P^1_k\setminus U$ is a union of $k$-points,
Hypothesis (HW) is expected to hold for an arbitrary number field $k$,
thus the restriction $k=\Q$ could potentially be dropped.
}\erem

\brem{\rm 
When $P(t)$ is irreducible, to prove Theorem \ref{adrasan2},
one can use the older fibration results \cite[Theorems 3.1, 4.1]{CSS98b} 
in place of Proposition \ref{h-w}.
}\erem

\section{Hasse principle for Kummer surfaces} \label{S4}

 \subsection{$2$-Selmer group in quadratic twist families} \label{ssec:twisting}
Let $k$ be a number field. Let $b\in k$ be such that $P(b)\neq 0$. We choose a finite set $S_{0,b}$   of places of $k$ containing $S_0$, all primes where $b$ is not integral and all primes dividing the discriminant $d_b=P(b)^2\,\discr(P(t))$ of $(t-b)P(t)$.  

Let $v$ be a place of $k$. To simplify notation,
we adopt the convention that the local invariant $\inv_v$ is a map $\Br(k_v)[2]\to \F_2$.

Recall from Section \ref{new geom} that $G$ is canonically isomorphic to $J_b[2]$. Composition of the Weil pairing and the local invariant map induces a non-degenerate pairing 
\begin{equation}\label{eq:local_tate_pairing}
\H^1(k_v,G)\times \H^1(k_v,G) \stackrel{\cup}{\longrightarrow}\textup{Br}(k_v)[2] \stackrel{\textup{inv}_v}{\longrightarrow}\F_2,
\end{equation}
which we refer to as the local Tate pairing. 
 Write $\del:J_b(k_v)\to \H^1(k_v,G)$ for the Kummer map. It is well known that
$\del(J_b(k_v))$ is its own orthogonal complement under the pairing \eqref{eq:local_tate_pairing},
see, e.g., \cite[Proposition 4.11]{PR21}.  

If  $v$ is an odd place, then $\dim_{\mathbb{F}_2} \del(J_b(k_v))=\dim_{\mathbb{F}_2} J_b(k_v)[2]$ by \cite[Lemma I.3.3]{Mil06}.
Thus $\dim_{\mathbb{F}_2}\H^1(k_v,G)= 2\,\dim_{\F_2} J_b(k_v)[2].$ 
Provided $v\notin S_0$, evaluation of cocycles at the Frobenius element 
$\textup{Frob}_v\in \textup{Gal}(k_v^{\textup{nr}}/k_v)$ induces an isomorphism
\begin{equation} \label{delta_frob_unram}
\H^1_\nr(k_v,G) \cong G/(\textup{Frob}_v-1).
\end{equation}

Let $F/k$ be a quadratic (or trivial) extension of $k$. For each place $v$ of $k$, we  view the Kummer map for $J_b^F$ as taking values in $\H^1(k_v,G)$ too, and denote by $\del^F:J_b^F(k_v)\to \H^1(k_v,G)$ the resulting map. The Weil pairing on $J_b[2]$ is unchanged under quadratic twist, so  $\del^F(J_b^F(k_v))$ is its own orthogonal complement with respect to the local Tate pairing also. 

Suppose that $v\notin S_{0,b}$. 
 If $F/k$ is unramified at $v$, then $J_b^F$ has good reduction at $v$, hence
\begin{equation}
\del^F(J_b^F(k_v))=\H^1_{\textup{nr}}(k_v,G). \label{bonn1}
\end{equation} 
On the other hand, if $v$ is  ramified in $F/k$, then by \cite[Lemma 4.3]{HS16} we have
\begin{equation}
\del^F(J_b^F(k_v))\cap \H^1_{\textup{nr}}(k_v,G)=0. \label{bonn2}
\end{equation}

We now turn to the global picture. The $2$-Selmer group $\Sel_2(J_b^F)$ is  the subgroup of 
$\H^1(k,G)$ consisting of elements $x$ such that $\loc_v(x)\in \del^F(J_b^F(k_v))$ for all places $v$ of $k$. It sits in a short exact sequence 
 \begin{equation} \label{eq:selmer_sequence}
 0\to J_b^F(k)/2 \to \Sel_2(J_b^F) \to \Sha(J_b^F)[2]\to 0,
 \end{equation}
 where $\Sha(J_b^F)$ is the Shafarevich--Tate group. We denote by 
\[\left \langle~,~\right \rangle_F:\Sha(J_b^F)\times \Sha(J_b^F) \to \mathbb{Q}/\mathbb{Z}\]
 the Cassels--Tate pairing. Since $C_b^F$ has a rational Weierstrass point, $\left \langle~,~\right \rangle_F$ is alternating by \cite[Theorem 8 and Corollary 12]{PS99}. If the $2$-primary subgroup of $\Sha(J_b^F)$ is finite then the same result shows that $\dim_{\mathbb{F}_2}\Sha(J_b^F)[2]$ is even.
   
 For a place $v$ of $k$, we denote by $\Sel_2(J_b^F)^{(v)}$ the subgroup of $\H^1(k,G)$ consisting  of elements $x$ such that $\loc_{w}(x)$ lies in $\del^F(J_b^F(k_{w}))$ for all places $w\neq v$.  That is, we impose the usual Selmer conditions at all places other than $v$ but impose no condition at $v$. The following result has its origin in work of Mazur--Rubin \cite{MR10}.

\ble \label{prop:variation_of_selmer_structure}
Let $F/k$ be a quadratic extension and  $v$  an odd place of $k$.

{\rm (i)}  The image of $\loc_v\big(\Sel_2(J_b^F)^{(v)}\big)$ in $\H^1(k_v,G)/\del^F(J_b^F(k_v))$ is the orthogonal complement of $\loc_v\big(\Sel_2(J_b^F)\big)$ with respect to the non-degenerate pairing 
\[\frac{\H^1(k_v,G)}{\del^F(J_b^F(k_v))} \times \del^F(J_b^F(k_v)) \to \F_2\]
induced by \eqref{eq:local_tate_pairing}. 

{\rm (ii)} Suppose that $v\notin\Ram(F/k)$ and $J_b^F$ has good reduction at $v$. Let $F'/k$ be a quadratic extension ramified at $v$ such that $\del^F(J_b^F(k_w))=\del^{F'}(J_b^{F'}(k_w))$ for every place $w\neq v$. Then
\[
\dim_{\mathbb{F}_2} \Sel_2(J_b^{F'})-\dim_{\mathbb{F}_2} \Sel_2(J_b^F)  =n_1-n_2, \]
 where $n_1= \dim_{\mathbb{F}_2}\loc_v\big(\Sel_2(J_b^{F'})\big)$ and $n_2=\dim_{\mathbb{F}_2}\loc_v\big(\Sel_2(J_b^F)\big)$. Moreover, writing $r= \dim_{\mathbb{F}_2}\big(G/(\textup{Frob}_{v}-1)\big)$, we have \[n_1+n_2\leq r\quad \textup{ and }\quad  n_1+n_2\equiv r \pmod 2.\]
\ele
{\em Proof.} 
(i) This is an immediate consequence of the Poitou--Tate duality.     

\noindent (ii) The first equality is clear. We have
$\del^F(J_b^F(k_v))\cap \del^{F'}(J_b^{F'}(k_v))=0$, as follows from 
(\ref{bonn1}) and (\ref{bonn2}). In this situation, using (\ref{delta_frob_unram}),
\cite[Lemma 3.27]{Har19} with $T=\{v\}$ gives both the inequality and the congruence.
 \hfill $\Box$
 
\medskip

The following lemma is a variant of \cite[Theorem 10.12]{Mor19} and will be used to control the parity of the dimension of the $2$-Selmer group in the family of quadratic twists of $J_b$.

\ble \label{lem:selmer_parity}
 Let $F=k(\sqrt{a})/k$ be a quadratic extension. For $v\in S_{0,b}$,   denote by $F_v$ the completion of $F$ at any place extending $v$ (the choice of which is not important).  
Then we have
\begin{equation} \label{selmer_parity}
(-1)^{\dim_{\mathbb{F}_2} \Sel_2(J_b^F)}=(-1)^{\dim_{\mathbb{F}_2} \Sel_2(J_b)}\cdot 
\prod_{v\in {S_{0,b}}}\left(\big(d_b, a\big)_{k_v}\cdot 
(-1)^{\dim_{\mathbb{F}_2}J_b(k_v)/N(J_b( F_v))}\right),
\end{equation}
where $N\colon J_b(F_v)\to J_b(k_v)$ is the norm map and $(~,~)_{k_v}$ is the
Hilbert symbol. 
\ele
{\em Proof.} 
If $v\notin S_{0,b}$ is a prime of $k$, then by \cite[Theorem 1.8]{Mor23a} the quantity 
 \[(d_b,a)_{k_v}\cdot (-1)^{\dim_{\mathbb{F}_2}J_b(k_v)/N(J_b( F_v))}\]
 is equal to the local root number of $J_b$ over $F_v$, which is $1$  since $J_b$ has good reduction at $v$ (since $C_b$ and $C^F_b$ have a rational Weierstrass point, the terms  $\epsilon(C_b/k_v)$ and $\epsilon(C^F_b/k_v)$ appearing in \cite[Conjecture 1.7]{Mor23a} are equal to $0$ for all places $v$). Consequently, we may replace the product over $v\in S_{0,b}$ in \eqref{selmer_parity} with the product over all places $v$ of $k$. The result now follows from \cite[Theorem 10.12]{Mor19} and the product formula for the Hilbert symbol. 
  \hfill $\Box$

 \subsection{The case when $P(t)$ has Galois group ${\rm S}_5$, $\AA_5$ or $\FF_{20}$} \label{SAF}

In this section we prove  Theorem \ref{thm} (a). The proof starts with the following key proposition.  

 \bpr \label{prop:carefully_chosen_twist}
Under the assumptions of Theorem \ref{thm} (a), if the Kummer surface $\Kum(J_{b,\delta})$ is everywhere locally soluble, then there is a quadratic extension $F/k$ such that  $J_{b,\delta}^{F}(\bA_k)\neq \emptyset$ and  $\dim_{\mathbb{F}_2} \Sel_2(J_b^F)$ is odd. 
 \epr
 {\em Proof.} 
Take a prime $w$ outside of $S_0$ such that $\val_w(P(b))=1$
We first show that if $F=k(\sqrt{a})$ is a quadratic extension in which $w$ is inert, then
$w$ contributes non-trivially to the right hand side of \eqref{selmer_parity}, that is,
\begin{equation} \label{local_root_number}
\big(d_b, a\big)_{k_w} \cdot(-1)^{\dim_{\mathbb{F}_2}J_b(k_w)/N(J_b(F_w))}=-1. 
\end{equation}
Indeed, up to squares, both $d_b$ and $a$ are units at $w$, hence $(d_b, a\big)_{k_w}=1$.
Next, our assumptions imply, as in the proof Proposition \ref{jan1}, that the group scheme of connected
components of the special fibre of the N\'eron model of $J_b$ is the constant group $\Z/2$.
Since $F$ is unramified at $w$, by \cite[Proposition I.3.8]{Mil06}, we have 
$\H^1_\nr(k_w,J^F_b)\cong\H^1(\F_w,\pi_0)\cong\Z/2$. From the exact sequence
$$0\to J_b^F\to R_{F/k}(J_b)\to J_b\to 0$$
we see that the cokernel of the norm map $J_b(F_w)\to J_b(k_w)$ is the kernel of 
the restriction map $\H^1(k_w,J_b^F)\to \H^1(F_w, J_b)$. Since $F$ is unramified at $w$,
this kernel is equal to the kernel of $\H^1_\nr(k_w,J^F_b)\to \H^1_\nr(F_w,J^F_b)$,
which is identified with multiplication by
2 on $\Z/2$, because $w$ is inert in $F$. Hence $J_b(k_w)/N(J_b(F_w))\cong\Z/2$, proving our claim.

We now adapt  \cite[Proposition 6.2]{HS16} (which has its origin in \cite[Section 5]{SSD05}) to construct a desired quadratic extension $F/k$. For convenience, we recall the setup of that result. Denote by $\tilde{J}_{b,\delta}$ the blow-up of $J_{b,\delta}$ in the $16$ fixed points of its antipodal involution $\iota$. Let $\mathcal{Y}$ denote the quotient of $\tilde{J}_{b,\delta}\times \mathbb{G}_{m,k}$ by the diagonal action of $\mu_2$ (acting by the extension of $\iota$ on the first factor, and by $-1$ on the second factor). The fibre of $\mathcal{Y}$ over $a\in \mathbb{G}_m(k)$ identifies with the quadratic twist of $\tilde{J}_{b,\delta}$ by $k(\sqrt{a})$. Let $\mathcal{Y}\subseteq \mathcal{X}$ be a smooth compactification fitting into a commutative diagram
$$\xymatrix{ \mathcal{Y}\ar[r]\ar[d]& \mathcal{X}
 \ar[d] \\ 
 \mathbb{G}_{m,k}\ar[r] & \mathbb{P}^1_k}$$
Since $\textup{Gal}(P)$ is ${\rm S}_5$, $\AA_5$ or $\FF_{20}$, it follows from  Lemma \ref{lemma:assumptions_a_and_b}  that $J_b$ satisfies conditions (a) and (b) of \cite[Section 2]{HS16}. Consequently, \cite[Proposition 6.1]{HS16} applies to show that the vertical Brauer group of $\mathcal{X}$ over $\mathbb{P}^1_k$ is the image of $\textup{Br}(k)$ in $\textup{Br}(\mathcal{X})$. 

By Proposition \ref{jan1}
local solubility of $\Kum(J_{b,\delta})$ at each place $v\in S_{0,b}\setminus \{w\}$ ensures the existence of a local extension $F_v=k_v(\sqrt{a_v})/k_v$ of degree at most $2$ such that $\tilde{J}^{F_v}_{b,\delta}$ has a $k_v$-point, say $P_v$. Set 
\[\kappa=(-1)^{\dim_{\mathbb{F}_2} \Sel_2(J_b)}\cdot \prod_{v\in S_{0,b}\setminus\{w\}}\left(\big(d_b, a_v\big)_{k_v}\cdot (-1)^{\dim_{\mathbb{F}_2}J_b(k_v)/N(J_b( F_v))}\right).\]
If $\kappa=-1$ then take $a_w=1$. If $\kappa=1$ then take $a_w\in k_w^\times$ 
to be any element generating the unique degree $2$ unramified extension of $k_w$. Let $F_w=k_w(\sqrt{a_w})$. Since $\Kum(J_{b,\delta})(k_w)\neq\emptyset$ and $F_w/k_w$ is unramified,
we have $\tilde{J}^{F_w}_{b,\delta}(k_w)\neq \emptyset$ by Proposition \ref{jan1}. Fix $P_w \in \tilde{J}^{F_w}_{b,\delta}(k_w)$.   Using local solubility of $\Kum(J_{b,\delta})$ at primes outside $S_{0,b}$, we extend the collection of local points $(P_v,\sqrt{a_v})_{v\in S_{0,b}}$ to an adelic point $(y_v)\in \mathcal{X}(\bA_k)$. 

Since the quotient of the vertical Brauer group of $\mathcal{X} \to \mathbb{P}^1_k$ by the image of $\textup{Br}(k)$ is trivial, we can use the fibration method to produce an adelic point $(y'_v)\in \mathcal{X}(\bA_k)$, arbitrarily close to $(y_v)$,  and such that the image of $(y_v')$ in $\mathbb{P}^1(\bA_k)$ is a $k$-point $a\in k^\times=\G_m(k)$. Specifically, we apply \cite[Theorem 9.17]{HW16} with $B=0$ and $U=\mathbb{G}_{m,k}$, noting that the cited result applies unconditionally thanks to \cite[Theorem 9.11]{HW16}.

Let $F=k(\sqrt{a})$. Then $J_{b,\delta}^F$ is everywhere locally soluble.  Further, for each $v\in S_{0,b}$, we can ensure that $a$ and $a_v$ are close enough that their square roots generate the same extension of $k_v$.  Having done this, if $\kappa=-1$ then $w$ splits in $F$ by construction and we see from Lemma \ref{lem:selmer_parity} that $\dim_{\mathbb{F}_2} \Sel_2(J_b^F)$ is odd. 
If $\kappa=1$ then $w$ is inert in $F$ and combining Lemma \ref{lem:selmer_parity} with \eqref{local_root_number} gives the same conclusion.
 \hfill $\Box$
 
 \medskip

 \noindent {\em Proof of Theorem \ref{thm} (a).} 
 We can assume without loss of generality that $\delta \neq 0$, since otherwise $\Kum(J_{b,\delta})$ contains a $k$-line.
 
 Suppose that $\Kum(J_{b,\delta})$ is everywhere locally soluble and let $F=k(\sqrt{a_0})$ be as in the statement of  Proposition \ref{prop:carefully_chosen_twist}. 
Then $\dim_{\mathbb{F}_2}\Sel_2(J_b^F)$ is odd and $\delta\in\Sel_2(J_b^F)$ since $J_{b,\delta}^F$ is everywhere locally soluble. 

If $\dim_{\mathbb{F}_2} \Sel_2(J_b^F)=1$ then we are done. Indeed, as in Section \ref{ssec:twisting}, finiteness of the $2$-primary part of $\Sha(J_b^F)$ implies that the $\mathbb{F}_2$-dimension of $\Sha(J_b^F)[2]$ is even. From the exact sequence \eqref{eq:selmer_sequence} we  conclude that $\delta$ maps to $0$ in $\Sha(J_b^F)$, hence the $2$-covering $J_{b,\delta}^F$  has a $k$-point. Thus so does $\textup{Kum}(J^F_{b,\delta})\cong  \textup{Kum}(J_{b,\delta})$.

Suppose  $\dim_{\mathbb{F}_2} \Sel_2(J_b^F)\geq 3$. We will construct another quadratic extension $F'/k$ such that $\delta \in \textup{Sel}(J_b^{F'})$ and 
\begin{equation}
\dim_{\mathbb{F}_2} \Sel_2(J_b^{F'})=\dim_{\mathbb{F}_2} \Sel_2(J_b^F)-2.
\label{bonn3}
\end{equation}
Let $T=\{\delta=\alpha_0,\alpha_1,...,\alpha_t\}$ be an $\mathbb{F}_2$-basis 
of $\Sel_2(J_b^F)$. The resulting map 
$$\varphi_T:\Ga\to G(\bar k)^{t+1}\rtimes \textup{Gal}(P)$$ is 
surjective by Lemmas \ref{lemma:assumptions_a_and_b} and  \ref{is_iso_vart},
where $q=2$, so that
$$\Gal(K_T/k)=G(\bar k)^{t+1}\rtimes \textup{Gal}(P).$$
   Fix a finite set $S$ of places of $k$ containing $S_{0,b}\cup\Ram(F/k)$. Fix 
$\sigma_0\in\textup{Gal}(P)$ acting on the roots of $P(t)$ as a double transposition, noting that
$\dim_{\mathbb{F}_2}G/(\sigma_0-1)=2$.
Further, fix elements $x_1$ and $x_2$ of $G$ whose images in $G/(\sigma_0-1)$ are 
$\mathbb{F}_2$-linearly independent, and let 
$$\sigma:=((0,0,...,x_1,x_2),\sigma_0)\in \textup{Gal}(k_T/k).$$ 
Our assumption that $\textup{Gal}(P)$ is 
${\rm S}_5$, $\AA_5$, or $\FF_{20}$ implies that every homomorphism $\textup{Gal}(P)\to \mathbb{Z}/2$ vanishes on $\si_0$. Thus by Lemma \ref{lem:existence_of_extensions} we can find an odd place $v\notin S$ such that $\Frob_v\in\textup{Gal}(k_T/k)$ lies in the conjugacy class of $\sigma$, and   $a\in k^\times$ such that $k(\sqrt{a})/k$ is ramified at $v$, split at all places in $S$, and unramified elsewhere. Recall that $F=k(\sqrt{a_0})$. Set $F'=k(\sqrt{aa_0})$. 
By construction, we have $\del^F(J_b^F(k_w))=\del^{F'}(J_b^{F'}(k_w))$ for every place $w\neq v$.  Our choice of $\sigma$ implies both that   $\loc_{v}(\delta)=0$, hence $\delta \in  \Sel_2(J_b^{F'})$, and that the image of $\Sel_2(J_b^F)$ in 
\[\H^1_{\textup{nr}}(k_v,G)\cong G/(\textup{Frob}_{v}-1)\]
is $2$-dimensional. Now Lemma \ref{prop:variation_of_selmer_structure} (ii) gives (\ref{bonn3}).
We now relabel $F'$ as $F$ and iterate the above argument, eventually arriving at a quadratic extension $F'/k$ such that $ \Sel_2(J_b^{F'})$  is generated by $\delta$. We can then conclude as before.
 \hfill $\Box$

\subsection{The Cassels--Tate pairing under quadratic twist}

In this section we study the variation of the Cassels--Tate pairing on the $2$-Selmer group
under quadratic twist, but first we need some more algebraic prerequisites. 

Let $M$ be a finite abelian group 
with a continuous action of $\Ga=\Gal(\bar k/k)$. Let us write $\Gal(M)$ for the image of 
$\Ga$ in $\Aut(M)$ and let $\varphi\colon \Ga\to \Gal(M)$ be the natural surjective map.
Let $T=\{\alpha_1,\ldots,\alpha_t\}\subset\H^1(k,M)$. 
For each $i=1,\ldots, t$ choose a 1-cocycle 
$\tilde\alpha_i\colon\Ga\to M$ representing $\alpha_i$.
As in Section \ref{1.1}, the map  
$$\varphi_T\colon\Ga\to M^T\rtimes\Gal(M), \quad \varphi_T(g)=\big(\tilde\alpha_1(g),\ldots,\tilde\alpha_t(g),\varphi(g)\big),$$
is a homomorphism. 
Let $k_T\subset\bar k$ be the fixed field of $\Ker(\varphi_T)$.
The field $k_T$  contains the smallest Galois extension of $k$ through which each of the cocycles $\tilde\alpha_1,\ldots,\tilde\alpha_t$ factors, thus
$T\subset\H^1(\Gal(k_T/k),M)\subset \H^1(k,M)$.

Suppose further that $M$ is a finite-dimensional $\F_2$-vector space
with a $\Ga$-equivariant non-degenerate bilinear paring 
$\cup\colon M\times M\to \F_2$. It induces a pairing
$$\cup\colon \H^1(\Gal(k_T/k),M)\times \H^1(\Gal(k_T/k),M)\to \H^2(\Gal(k_T/k),\F_2).$$
In particular, 
the cohomology class $\alpha_i\cup\alpha_j$ is given by the normalised 2-cocycle
$\tilde \alpha_i\cup\tilde\alpha_j\in Z^2(\Gal(k_T/k),\F_2)$ that
sends $\si,\tau\in \Gal(k_T/k)$ to $\tilde{\alpha}_i(\sigma)\cup\sigma\tilde{\alpha}_j(\tau)$.

Let $\F_2^T$ be an $\F_2$-vector space of dimension $t=|T|$ with basis $e_1,\ldots, e_t$.
Let $\tilde\xi_T\in Z^2(\Gal(k_T/k),\wedge^2(\F_2^T))$ be the normalised 2-cocycle
whose $e_i\wedge e_j$-coordinate is 
$\tilde \alpha_i\cup\tilde\alpha_j\in Z^2(\Gal(k_T/k),\F_2)$.
Let $\xi_T\in \H^2(\Gal(k_T/k),\wedge^2(\F_2^T))$ be the class of $\tilde\xi_T$.

Consider the associated central extension
$$0\to \wedge^2(\F_2^T)\to E_T\to \Gal(k_T/k)\to 0.\eqno{(E_T)}$$
The underlying set of $E_T$ is $\wedge^2(\F_2^T)\times \Gal(k_T/k)$ 
with the group operation defined so that 
the $e_i\wedge e_j$-coordinate of the product of $\si,\tau\in \Gal(k_T/k)$ is
$\tilde{\alpha}_i(\sigma)\cup\sigma\tilde{\alpha}_j(\tau)$.

The next statement generalises \cite[Corollary 4.18]{Mor23b}. 
We give a self-contained proof inspired by the proof of \cite[Proposition 4.11]{Har19}.

\bpr \label{rosenmontag}
Assume that $\xi_T$ goes to zero under the inflation map 
$$\H^2(\Gal(k_T/k),\wedge^2(\F_2^T))\to\H^2(k,\wedge^2(\F_2^T)),$$ 
so that there is a normalised $1$-cochain
$\gamma\in C^1(k,\wedge^2(\F_2^T))$, $\gamma(e)=0$, 
such that $\tilde\xi_T=d\gamma$. Then the function
$$\hat\varphi_T\colon \Ga\to E_T, \quad \quad 
\hat\varphi_T(g)=(\gamma(g),\varphi_T(g)),$$
is a homomorphism. 
If $\varphi_T(\Ga)=M^T\rtimes\Gal(M)$, then $\hat\varphi_T(\Ga)=E_T$
and we have a natural isomorphism of abelianisations 
$(E_T)^{\rm ab}\cong \Gal(k_T/k)^{\rm ab}$.
\epr
{\em Proof.} The pullback of $(E_T)$ under the sujective map
$\Ga\to\Gal(k_T/k)$ given by $\varphi_T$ is split, hence this map lifts to a homomorphism
$\Ga\to E_T$. The choice of this lifting is equivalent to the choice of a normalised $1$-cochain 
$\gamma$ such that $\tilde\xi_T=d\gamma$, 
so that the lifting $\Ga\to E_T$ is the homomorphism $(\gamma(g),\varphi_T(g))$.

Assume that $\varphi_T\colon\Ga\to M^T\rtimes\Gal(M)$ is surjective.
To prove that $\hat\varphi_T(\Ga)=E_T$ it is enough to prove that
the restriction of $\hat\varphi_T$ to $\Ga_{k_T}=\Gal(\bar k/k_T)$ is a
surjective homomorphism $\Ga_{k_T}\to\wedge^2(\F_2^T)$. This homomorphism sends
$g\in \Ga_{k_T}$ to $\gamma(g)$.
It is enough to produce for any $i\neq j$ an element $g\in \Ga_{k_T}$
such that $\gamma_{i,j}(g)=1$ and $\gamma_{k,l}(g)=0$ when $(k,l)\neq (i,j)$.

Take $x,y\in M$ such that $x\cup y=1$. 
Since $\varphi_T\colon\Ga\to M^T\rtimes\Gal(M)$ is surjective, we can find $\si\in\Ga$
such that $\varphi_T(\si)$ is the element of the subgroup $M^T$
of $M^T\rtimes\Gal(M)$ whose $i$-th coordinate is $x$ and all other coordinates are $0$.
Likewise, we can find $\tau\in\Ga$
such that $\varphi_T(\tau)\in M^T\subset M^T\rtimes\Gal(M)$ has $y$ in the $j$-th coordinate, and $0$ in all other coordinates. Since $M^T$ is abelian,
we obtain $\si\tau\si^{-1}\tau^{-1}\in \Ga_{k_T}$. 
Note that $\si, \tau$ are in $\Ga_{M}:=\Ker[\Ga\to \Gal(M)]$ which acts trivially on $M$.
Thus each $\tilde\alpha_i$ is a homomorphism $\Ga_M\to M$. Likewise,
since $\tilde\xi_T=d\gamma$, for any $g,h\in \Ga_M$ we have
$$\tilde\alpha_i(g)\cup\tilde\alpha_j(h)=\gamma(gh)-\gamma(g)-\gamma(h).$$
Using that $\tilde\alpha_i(\si)=x$, $\tilde\alpha_j(\si)=0$, $\tilde\alpha_i(\tau)=0$, 
$\tilde\alpha_j(\tau)=y$, together with $\gamma(e)=0$, we compute
$\gamma_{i,j}(\si\tau\si^{-1}\tau^{-1})=x\cup y=1$. A similar calculation gives 
$\gamma_{k,l}(\si\tau\si^{-1}\tau^{-1})=0$ for $(k,l)\neq (i,j)$.

As a by-product, we obtain that the subgroup $\wedge^2(\F_2^T)\subset E_T$ 
is contained in the commutator subgroup $[E_T,E_T]$. 
In particular, the abelianisation of $E_T$ identifies with the abelianisation of $\Gal(k_T/k)$.
\hfill $\Box$

\medskip

We now turn to the variation of the Cassels--Tate pairing on the $2$-Selmer group
under quadratic twist. 
With a view to potential applications we work with an arbitrary principally polarised abelian variety
$A$ over a number field $k$.

The Weil pairing $e_2\colon A[2]\times A[2] \to \F_2$ associated to the given principal polarisation
is $\Ga$-equivariant, non-degenerate, and alternating. We denote by
$$\cup_{e_2}\colon C^i(k,A[2])\times C^j(k,A[2])\to C^{i+j}(k,\F_2)$$ 
the induced cup-product on cochains or cohomology classes. In particular,
we have a pairing
$$\cup_{e_2}\colon \H^1(k,A[2])\times \H^1(k,A[2])\to \H^2(k,\F_2)$$
compatible with a similar pairing for $k_v$ under the restriction maps.

Let $\alpha,\beta\in\Sel_2(A)\subset\H^1(k,A[2])$.
We apply the construction from the beginning of this section to $M=A[2]$ and 
$T=\{\alpha,\beta\}$ and so obtain an extension 
$$0\to\F_2\to E_T\to \Gal(k_T/k)\to 0.$$
Recall that for each place $v$ of $k$, 
the subgroup $\partial(A(k_v))\subset\H^1(k_v,A[2])$ is isotropic with respect to the local Tate pairing. Using this together with global reciprocity for the Brauer group
we see that $\alpha \cup_{e_2} \beta=0$. Choosing 2-cocycles 
$\tilde\alpha,\tilde\beta\in Z^2(k,A[2])$ representing $\alpha,\beta$, respectively,
we can thus find a 1-cochain $\gamma\in C^1(k,A[2])$ such that 
$d\gamma=\tilde\alpha\cup_{e_2}\tilde\beta$. As in Proposition \ref{rosenmontag}
we obtain a homomorphism $\hat\varphi_T\colon \Ga\to E_T$ given by
$(\gamma,\varphi_T)$. Let $K_T$ be the fixed field of $\Ker(\hat\varphi_T)$.

Let $\left \langle~,~\right \rangle_2^{\textup{CT}}$ be the $\F_2$-valued pairing on 
$\textup{Sel}_2(A)$ given by pulling back the Cassels--Tate pairing along 
the surjective map $\textup{Sel}_2(A)\to \Sha(A)[2]$. For $F/k$ quadratic, we denote by $\left \langle ~,~\right \rangle_{2,F}^{\textup{CT}}$ the corresponding pairing on $\textup{Sel}_2(A^F)$. The sum of these pairings gives a pairing on $\textup{Sel}_2(A)\cap \textup{Sel}_2(A^F)$. The following result describing this pairing in favourable circumstances is a direct generalisation of \cite[Proposition 3.3]{Mor23b} which was inspired by \cite[Proposition 3.29]{Har19}
and \cite[Theorem 3.2]{Smi16}).


Recall that for a finite field extension $k'/k$ we write $\Ram(k'/k)$ for the set of primes of $k$ ramified in $k'$.
 
\bpr \label{prop:variation_of_CTP} 
Let $A$ be a principally polarised abelian variety over a number field $k$. Let
$T=\{\alpha, \beta\} \subset\Sel_2(A)$. 
Let $S$ be a finite set of places containing all archimedean places, all places dividing $2$, 
all places of bad reduction for $A$, and $\Ram(K_T/k)$, where $K_T$ is 
the finite Galois extension of $k$ constructed above.
Let $F/k$ be a quadratic extension that splits at all places of $S$ and such that
$\loc_{v}(\alpha)=\loc_{v}(\beta)=0$ for all $v\in \textup{Ram}(F/k)$. 
For all $v\in \Ram(F/k)$, choose 
$P_v, Q_v\in A[2](\ov k_v)$ such that $dP_v=\loc_{v}(\tilde{\alpha})$ and 
$dQ_v=\loc_{v}(\tilde{\beta})$. 
 Then $\alpha, \beta \in \Sel_2(A^F)$ and we have 
 \begin{equation} \label{eq:CT_variation}
 \left\langle \alpha,\beta \right \rangle_{2}^{\textup{CT}}+\left \langle \alpha,\beta \right\rangle_{2,F}^{\textup{CT}}=\sum_{v\in \textup{Ram}(F/k)} \big(e_2(P_v,(\textup{Frob}_v-1)Q_v)-\tilde{\gamma}(\textup{Frob}_v)\big) .
 \end{equation}
 \epr
{\em Proof.} 
Our assumption that $F/k$ splits at all places of $S$ ensures that the $2$-Selmer conditions for $A$ and $A^F$ agree at all places $v\notin \textup{Ram}(F/k)$. Since 
$\loc_{v}(\alpha)=\loc_{v}(\beta)=0$ at all places $v\in \textup{Ram}(F/k)$, it follows that $\alpha,\beta \in \textup{Sel}_2(A^F)$.

\smallskip

There is a pairing on $\Sel_2(A)\cap \Sel_2(A^F)$ constructed in
\cite[Definition 5.6]{Mor19} and shown in \cite[Lemma 5.8]{Mor19} to be equal to the sum of the Cassels--Tate pairings for $A$ and $A^F$. We now recall \cite[Definition 5.6]{Mor19}. 
 
Write $\chi:\Ga\to \mu_2$ for the quadratic character of $F/k$. We write $\chi_v:=\loc_v(\chi)$.

We denote the unique isomorphism $\mu_2\to \F_2$ by $[x]$, so that $[1]=0$ and $[-1]=1$.

The necessary global choices of $1$-cocycles $\tilde\alpha, \tilde\beta$ and 
a $1$-cochain $\gamma$ have been already done, our pairing
is a sum of local terms, one for each place $v$ of $k$. 

Since $\alpha \in \textup{Sel}_2(A)$, we can find  a point $P_v\in A(\ov k_v)$ such that
$dP_v=\loc_{v}(\tilde{\alpha})$. When $v\in \textup{Ram}(F/k)$,
we have $\loc_v(\alpha)=0$; for such places $v$ we take
$P_v\in A[2](\ov k_v)$ as in the proposition.
Pick $P_v'\in A(\ov k_v)$ such that $2P_v'=P_v$.  
Then $\rho_v:=dP'_v$ is a $1$-coboundary so that 
for $\si\in\Gal(\ov k_v/k_v)$ we have $\rho_v(\si)=\si(P'_v)-P'_v$.

Doing the same for $A^F$ we find points $P_{v,\chi}, P'_{v,\chi}\in A(\ov k_v)$
such that $d_\chi P_{v,\chi}=\loc_{v}(\tilde{\alpha})$ and $P_{v,\chi}=2P'_{v,\chi}$.
When $v\in \textup{Ram}(F/k)$, we take $P_{v,\chi}=P_v$ and $P'_{v,\chi}=P'_v$.
The 1-coboundary $\rho_{v,\chi}:=d_\chi P'_{v,\chi}$ satisfies
$\rho_{v,\chi}(\si)=\chi_v(\si)\si(P'_{v,\chi})-P'_{v,\chi}$ for any $\si\in\Gal(\ov k_v/k_v)$.

Note that $\rho_v+\rho_{v,\chi}$ has values in $A[2](\ov k_v)$, so
$\rho_v+\rho_{v,\chi}\in C^1(k_v,A[2])$.
We claim that we have $d(\rho_v+\rho_{v,\chi})=[\chi_v]\cup \loc_v(\tilde{\alpha})$
in $Z^2(k_v,A[2])\subset Z^2(k_v,A)$. Since $2\rho_{v,\chi}= \loc_v(\tilde{\alpha})$ and
$d\rho_v=d(dP'_v)=0$, it is enough to check that $d\rho_{\chi,v}=[\chi_v]\cup (2\rho_{v,\chi})$.
Rewriting $d_\chi\rho_{\chi,v}=d_\chi(d_\chi P'_{v,\chi})=0$,  
we obtain that for any $\si,\tau\in \Ga$ we have 
$$0=\rho_{v,\chi}(\si)+\chi_v(\si)\si\rho_{v,\chi}(\tau)-\rho_{v,\chi}(\si\tau)=
d\rho_{v,\chi}(\si,\tau)-\si\rho_{v,\chi}(\tau)+\chi_v(\si)\si\rho_{v,\chi}(\tau),$$
which gives the desired equality $d\rho_{\chi,v}=[\chi_v]\cup (2\rho_{v,\chi})$, proving the claim.

The claim implies that 
$$(\rho_v+\rho_{v,\chi})\cup_{e_2}\loc_v(\tilde \beta)-[\chi_v]\cup\loc_v(\tilde{\gamma})$$ 
is a 2-cocycle with values in $\mu_2$. By \cite[Definition 5.6]{Mor19}
the local term of our pairing at $v$ is obtained by applying the local invariant $\inv_v$ to 
the cohomology class of this cocycle. 

\smallskip

Let us show that this local term at a place $v$ is zero unless $v\in\Ram(F/k)$.

If $v\in S$, then $[\chi_v]=0$. Thus $A\times_kk_v\cong A^F\times_kk_v$, so we can take
$P_v=P_{v,\chi}$ and $P'_v=P'_{v,\chi}$. Then $\rho_v+\rho_{v,\chi}=0$, so the local term 
of our pairing at $v$ is zero.

If $v\notin S\cup \Ram(F/k)$, then $\chi_v$ is unramified. Since $A$ has good reduction at $v$,
and $v$ does not divide 2, the classes $\alpha$ and $\beta$ are unramified at $v$. Thus 
without loss of generality we can assume that $\loc_v(\tilde\alpha)$ and $\loc_v(\tilde\beta)$
are 1-cocycles of $\Gal(k_v^\nr/k_v)$, and hence that $\loc_v(\gamma)$ is 
a 1-cochain $\Gal(k_v^\nr/k_v)\to A[2]$. This implies that $[\chi_v]\cup\loc_v(\tilde{\gamma})$
is a 2-cochain on $\Gal(k_v^\nr/k_v)\to A[2]$. Next, since $A(k_v^\nr)$ is 2-divisible, we can choose
points $P_v$ and $P'_v$ in $A(k_v^\nr)$, so that  $\rho_v$ is a 1-cochain on  $\Gal(k_v^\nr/k_v)$.
The same holds for $\rho_{v,\chi}$.
Moreover, the above calculation shows that we have an equality
$d(\rho_v+\rho_{v,\chi})=[\chi_v]\cup \loc_v(\tilde{\alpha})$ of 2-cochains on $\Gal(k_v^\nr/k_v)$.
Since $\H^2(k_v^\nr,\F_2)=0$, we conclude that the local term of our pairing at $v$ is zero.

Finally, let $v\in \Ram(F/k)$. In this case we have $P_v=P_{v,\chi}\in A[2](\ov k_v)$ and
$P'_{v,\chi}=P'_v\in A[4](\ov k_v)$. One immediately checks that $\rho_v-\rho_{v,\chi}=
[\chi_v]\cup P_v$. From this we obtain
$$(\rho_v+\rho_{v,\chi})\cup_{e_2} \loc_{v}(\tilde{\beta})+
\loc_v([\chi]\cup\tilde{\gamma})=
\big(([\chi_v] \cup P_v)-\loc_{v}(\tilde{\alpha})\big)\cup_{e_2} \loc_{v}(\tilde{\beta})-
\loc_v([\chi]\cup\tilde{\gamma}).$$
Both $\loc_{v}(\tilde{\alpha})$ and $\loc_{v}(\tilde{\beta})$ 
are $1$-cocycles representing elements of   $\partial(A(k_v))$. Since $\partial(A(k_v))$ is its own orthogonal complement under the local Tate pairing, we have $\textup{inv}_v(\loc_{v}(\tilde{\alpha}) \cup_{e_2} \loc_{v}(\tilde{\beta}))=0$. Thus the contribution  from the place $v$ is equal to 
$$\textup{inv}_v\big([\chi_v] \cup (P_v \cup_{e_2}\loc_{v}(\tilde{\beta}) -\loc_{v}(\tilde{\gamma}))\big).$$
This equals
$e_2(P_v,\tilde{\beta}(\Frob_v))-\tilde{\gamma}(\Frob_v)$,
since $P_v \cup_{e_2}\loc_{v}(\tilde{\beta}) -\loc_{v}(\tilde{\gamma})$ is an unramified quadratic character, while $\chi_v=\loc_{v}(\chi)$ is a ramified quadratic character.
 Noting that $\tilde{\beta}(\textup{Frob}_v)=(dQ_v)(\textup{Frob}_v)=(\textup{Frob}_v-1)Q_v$ and summing over $v\in \textup{Ram}(F/k)$ gives the result.   \hfill $\Box$

\subsection{The case when $P(t)$ has Galois group $\DD_{10}$} \label{D10}

In this section we prove Theorem \ref{thm} (b).
The proof proceeds along roughly the same lines as the proof of Theorem \ref{thm} (a). However, there are two complications. Firstly, as the action of the Galois group $\Ga$ on $G$ factors through
$\Gal(P)$, the fact that $\End_k(G)=\textup{End}_{\DD_{10}}(G)=\mathbb{F}_4$ 
implies that, given $\mathbb{F}_2$-linearly independent elements $\alpha_0,...,\alpha_t$ of $\H^1(k,G)$, it need no longer be the case that the smallest Galois extension through which all the $\alpha_i$ factor has Galois group $G(\bar k)^{t+1}\rtimes \textup{Gal}(P)$. By Lemma \ref{is_iso_vart}, one needs $\alpha_0,...,\alpha_t$ to be $\mathbb{F}_4$-linearly independent to guarantee this. Secondly, it is no longer the case that every homomorphism $\textup{Gal}(P)\to \mathbb{Z}/2\mathbb{Z}$ vanishes on double transpositions. The result is that, upon twisting, we are only able to guarantee that the $\mathbb{F}_2$-dimension of the  Selmer group decreases by $4$ at each step. We then need an additional argument to treat the case where the twisting process terminates in a $3$-dimensional Selmer group. To handle this, we introduce a second descent step. 
 
 We begin with two preparatory results which help handle the first difficulty. From now on we take all the notation from the statement of Theorem \ref{thm} (b). 
 
 \ble \label{lem:span_of_delta}
 Let $w\in \{w_1, w_2\}$ and suppose that $F/k$ is a quadratic extension unramified at $w$. If $\alpha\in \Sel_2(J_b^F)$ has $\loc_w(\alpha)\neq 0$, then 
$\mathbb{F}_4 \alpha\cap \Sel_2(J_b^F)=\{0,\alpha\}$.
 \ele
 
\noindent {\em Proof.}
By assumption, $\Frob_w$ acts on the roots of $P(t)$ as a double transposition. In particular,   we have
 \[ \dim_{\mathbb{F}_2}\del^F(J_b^F(k_w))=\dim_{\mathbb{F}_2}J_b^{F}(k_w)[2]=2. \]
The natural maps $\H^1(k,G)\to \H^1(k_w,G)\to \H^1(k_w^\nr, G)$ are 
linear transformations of $\F_4$-vector spaces. Suppose for contradiction that $\mathbb{F}_4\alpha\subset \Sel_2(J_b^F)$. 
Since $\loc_w(\alpha)\neq 0$, we then have 
$\del^F(J_b^F(k_w))=\mathbb{F}_4\loc_w(\alpha)$. As $\H^1_{\textup{nr}}(k_w,G)$ is an $\mathbb{F}_4$-vector space, 
$\del^F(J_b^F(k_w))\cap \H^1_{\textup{nr}}(k_w,G)$ is an $\F_4$-vector space too.
From Lemma \ref{mumbai3}, taking into account that $\pi_0\cong\Z/2$, we see that this intersection is $\Z/2$,  giving the sought contradiction. \hfill $\Box$
 
\ble \label{arranging_linearly_independent}
Let $F/k$ be a quadratic extension unramified at $w_1$ and $w_2$. Suppose that $\delta \in \Sel_2(J_b^F)$ and $\dim_{\mathbb{F}_2}\Sel_2(J_b^F)=3$. Then there is a quadratic extension $F'/k$  unramified at $w_1$ and $w_2$, such that $\delta \in \Sel_2(J_b^{F'})$ and such that either $\dim_{\mathbb{F}_2}\Sel_2(J_b^{F'})=1$, or $\dim_{\mathbb{F}_2}\Sel_2(J_b^{F'})=3$ and the vectors of every $\mathbb{F}_2$-basis of $\Sel_2(J_b^{F'})$ are $\mathbb{F}_4$-linearly independent in $\H^1(k,G)$.
\ele

\noindent {\em Proof.}
Suppose that $\Sel_2(J_b^F)$ has an $\mathbb{F}_2$-basis which is not $\mathbb{F}_4$-linearly independent in $\H^1(k,G)$. Let $V=\Sel_2(J_b^F)$
and let $W$ be the $\F_4$-linear span of $V$ in $\H^1(k,G)$. 
We have $\dim_{\F_2}(W)=2\dim_{\F_4}(W)\leq 4$.  
Let $x\in \mathbb{F}_4\setminus \mathbb{F}_2$.
As $\dim_{\F_2}(V)=3$ we see that $V\cap x V$ is non-zero and invariant under multiplication by $x$, so
it is a non-zero $\F_4$-vector space. Thus $\Sel_2(J_b^F)$ contains $\F_4\alpha_1$ for some $\alpha_1\neq 0$. 
Lemma \ref{lem:span_of_delta} implies that $\loc_{w_1}(\alpha_1)=0$ and 
$\loc_{w_2}(\alpha_1)=0$. By assumption
$\loc_{w_1}(\delta)\neq 0$ and $\loc_{w_2}(\delta)=0$, 
hence $\Sel_2(J_b^F)=\mathbb{F}_4\alpha_1 \oplus \mathbb{F}_2\delta$ and so
$\loc_{w_2}\big(\Sel_2(J_b^F)\big)=0$. 

By part (i) of Proposition \ref{prop:variation_of_selmer_structure}, a consequence of the vanishing of $\loc_{w_2}\big(\Sel_2(J_b^F)\big)$ is the existence of an element 
$\alpha_2 \in\Sel_2(J_b^F)^{(w_2)}$ with $\loc_{w_2}(\alpha_2)\neq 0$. 
The set $T=\{\delta, \alpha_1,\alpha_2\}$ is $\mathbb{F}_4$-linearly independent, thus by Lemma \ref{is_iso_vart} the associated map   
$\varphi_T:\textup{Gal}(k_T/k)\to G(\bar k)^3\rtimes\textup{Gal}(P)$ is an isomorphism.

Choose an $\F_4$-basis $\{x_1,x_2\}$ of $G$ and define
$\si=\varphi_T^{-1}\big(((0,x_1,x_2),e)\big) \in \Gal(k_T/k)$. 
Let  $S$ be a finite set of places of $k$ containing $S_{0,b}$ and $\Ram(F/k)$.  By Lemma \ref{lem:existence_of_extensions} we can find a place $v\notin S$ whose Frobenius element in $\textup{Gal}(k_T/k)$ lies in the conjugacy class of $\sigma$, and   
$a\in k^\times$ such that $k(\sqrt{a})/k$ is ramified at $v$, split at all places of $S$, 
and unramified elsewhere. 
Write $F=k(\sqrt{a_0})$ and set $F'=k(\sqrt{aa_0})$, noting that $F'/k$ is unramified at $w_1$ and $w_2$. Since  $\textup{Frob}_{v}$ has trivial image in $\textup{Gal}(P)$ by construction, evaluation of cocycles at $\textup{Frob}_{v}$ gives an isomorphism
$\H^1_{\textup{nr}}(k_v,G)\cong G.$
Our choice of $\sigma$ ensures that $\loc_{v}(\delta)=0$,
hence $\delta \in \Sel_2(J_b^{F'})$. It also ensures that
$\H^1_{\textup{nr}}(k_v,G)=
\F_4\loc_{v}(\alpha_1)\oplus\F_4\loc_{v}(\alpha_2)$. In particular, we have
$\dim_{\mathbb{F}_2}\loc_{v}\big(\Sel_2(J_b^{F})\big)=2$.

 By Lemma \ref{prop:variation_of_selmer_structure} (ii) we have
$\dim_{\F_2}\Sel_2(J_b^{F'})=1+n_1,$
where $n_1=\dim_{\mathbb{F}_2}\loc_{v}\big(\Sel_2(J_b^{F'})\big)$ is either $0$ or $2$.  If $n_1=0$  then we are done, so suppose that $n_1=2$. By construction, $\Sel_2(J_b^{F'})$ is contained in $\Sel_2(J_b^F)^{(v)}$, hence $\loc_{v}\big(\Sel_2(J_b^{F'})\big)$ and $\loc_{v}\big(\Sel_2(J_b^{F})\big)$ are orthogonal with respect to the local Tate pairing at $v$. Since $\H^1_{\textup{nr}}(k_v,G)$ is its own orthogonal complement with respect to this pairing, and since $\delta^{F'}(J_b^{F'}(k_{v}))$ and $\H^1_{\textup{nr}}(k_v,G)$ have trivial intersection,  the restriction of \eqref{eq:local_tate_pairing} to a pairing 
\begin{equation} \label{eq:restricted_tate_pairing}
\H^1_{\textup{nr}}(k_{v},G)\times \delta^{F'}(J_b^{F'}(k_v))\to \F_2
\end{equation}
is non-degenerate. From this we conclude that the $2$-dimensional $\mathbb{F}_2$-vector spaces $\loc_{v}\big(\Sel_2(J_b^{F})\big)$ and $\loc_{v}\big( \Sel_2(J_b^{F'})\big)$ are orthogonal complements with respect to the pairing \eqref{eq:restricted_tate_pairing}. Since $\loc_{v}(\alpha_2)\in \H^1_{\textup{nr}}(k_v,G)$   does not lie in $\loc_{v}\big(\Sel_2(J_b^{F})\big)$, we deduce the existence of an element $\eta\in\Sel_2(J_b^{F'})$ such that $\loc_{v}(\eta)$ pairs non-trivially with $\loc_{v}(\alpha_2)$. Since both $\eta$ and $\alpha_2$ satisfy the Selmer conditions for $J_b^F$ away from $w_2$ and $v$, we conclude from reciprocity for the Brauer group of $k$ that $\loc_{w_2}(\eta)$ pairs non-trivially with $\loc_{w_2}(\alpha_2)$ under the local Tate pairing at $w_2$. In particular, $\loc_{w_2}(\eta)\neq 0$. Thus 
$\Sel_2(J_b^{F'})\cong(\F_2)^3$ contains $\F_2\delta\oplus\F_2\eta$ and we have $\loc_{w_2}(\eta)\neq 0$. By the argument in the beginning of the proof, if $\Sel_2(J_b^{F'})$ has an 
$\mathbb{F}_2$-basis which is not $\mathbb{F}_4$-linearly independent in $\H^1(k,G)$, 
then we have $\loc_{w_2}\big(\Sel_2(J_b^{F'})\big)=0$, a contradiction. \hfill $\Box$

 \medskip
 We now give the analogue of  Proposition  \ref{prop:carefully_chosen_twist}.
 
 \bpr \label{prop:carefully_chosen_twist_2}
Under the assumptions of Theorem \ref{thm} (b), if the Kummer surface $\Kum(J_{b,\delta})$ is everywhere locally soluble, then there is a quadratic extension $F/k$, unramified at $w_1$ and $w_2$, such that  $J_{b,\delta}^{F}(\bA_k)\neq \emptyset$ and  $\dim_{\mathbb{F}_2} \Sel_2(J_b^F)$ is odd.  
 \epr
  \noindent {\em Proof.}
 We follow the proof of Proposition \ref{prop:carefully_chosen_twist} verbatim, with $w_2$  in place of $w$,  until it comes to showing that the vertical Brauer group of $\mathcal{X}$ over $\mathbb{P}^1_k$ is the image of $\textup{Br}(k)$ in $\textup{Br}(\mathcal{X})$. Since $\textup{Gal}(P)$ is now $\DD_{10}$, condition (b)  of \cite[Proposition 6.1]{HS16} is satisfied,  but condition (a) is not. 
 However, writing $T=\{\delta\}$, by Lemma \ref{is_iso_vart} it is nevertheless the case that the resulting map $\varphi_T:\textup{Gal}(k_T/k)\to G(\bar k)\rtimes \textup{Gal}(P)$ is an isomorphism, which is enough for the proof of \cite[Proposition 6.1]{HS16} to carry over unchanged to show that the vertical Brauer group of $\mathcal{X} \to \mathbb{P}^1_k$ is  the image of $\textup{Br}(k)$. We then continue to follow the proof of Proposition \ref{prop:carefully_chosen_twist} to construct the sought extension $F/k$, with the only additional input needed being the existence of an  {\em unramified} local extension $F_{w_1}/k_{w_1}$ such that  $\tilde{J}^{F_{w_1}}_{b,\delta}$ has a $k_{w_1}$-point. Such an extension is afforded by \cite[Lemma 3.16]{Har19}.
  \hfill $\Box$
  \bigskip
  
\noindent {\em Proof of Theorem \ref{thm} (b).} 
Suppose that $\Kum(J_{b,\delta})$ is everywhere locally soluble and let $F/k$ be as in the statement of Proposition \ref{prop:carefully_chosen_twist_2}, noting that $\delta \in   \Sel_2(J_b^F)$. Suppose that $\dim_{\mathbb{F}_2} \Sel_2(J_b^F)\geq 5$. We claim that we can find a further quadratic extension $F'/k$, unramified at $w_1$ and $w_2$, such that $\delta \in  \Sel_2(J_b^{F'})$ and 
\[\dim_{\mathbb{F}_2} \Sel_2(J_b^{F'})=\dim_{\mathbb{F}_2} \Sel_2(J_b^{F}) -4.\]
To begin, denote by $V$ the subgroup of  $\Sel_2(J_b^F)$ consisting of elements $\alpha$ for which $\mathbb{F}_4\alpha \subseteq \Sel_2(J_b^F)$. This is an $\mathbb{F}_4$-vector space. Let $\alpha_1,...,\alpha_n$ be an $\mathbb{F}_4$-basis for $V$. Lemma \ref{lem:span_of_delta} and the assumption that $\loc_{w_1}(\delta)\neq 0$ imply that $\delta \notin V$. Fix $x\in \mathbb{F}_4\setminus \mathbb{F}_2$, write $\alpha_0=\delta$, and extend $\alpha_0, \alpha_1, x\alpha_1,...,\alpha_n,x\alpha_n$ to an $\mathbb{F}_2$-basis for $\Sel_2(J_b^F)$, adding in elements  $\alpha_{n+1},...,\alpha_t$. The set $T=\{\alpha_0,\alpha_1,...,\alpha_t\}$ is then $\mathbb{F}_4$-linearly independent. By Lemma \ref{is_iso_vart}, the resulting map $\varphi_T:\textup{Gal}(k_T/k)\to G(\bar k)^{t+1}\rtimes \textup{Gal}(P)$ is an isomorphism. 
 Since $\dim_{\mathbb{F}_2} \Sel_2(J_b^F)\geq 5$, we can thus find $\sigma \in \textup{Gal}(k_T/k)$ with trivial image in $\textup{Gal}(P)$, such that $\alpha_0(\sigma)=0$ and  the elements
\[\alpha_1(\sigma), x\alpha_1(\sigma),...,\alpha_n(\sigma),x\alpha_n(\sigma), \alpha_{n+1}(\sigma),...,\alpha_{t}(\sigma)\]
span $G$. (Note that since $\sigma$ has trivial image in $\textup{Gal}(P)$  we can unambiguously evaluate cocycle classes at $\sigma$.) For example, if $n\geq 2$ then we can take $x_1,x_2$ to be an $\mathbb{F}_4$-basis for $G$ and take $\sigma=\varphi_T^{-1}\big(((0,x_1,x_2,0,...,0),e)\big)$. 

Fix a finite set $S$ of places of $k$ containing $S_{0,b}\cup\Ram(F/k)$.  By Lemma \ref{lem:existence_of_extensions} there is a place $v\notin S$ whose Frobenius element in $\textup{Gal}(k_T/k)$ lies in the conjugacy class of $\sigma$, and   $a\in k^\times$ such that $k(\sqrt{a})/k$ is ramified at $v$, split at all places in $S$, and is unramified otherwise.  Write $F=k(\sqrt{a_0})$ and set $F'=k(\sqrt{aa_0})$, noting that $F'/k$ is unramified at $w_1$ and $w_2$.  
Then $\loc_{v}(\delta)=0$ and $\dim_{\mathbb{F}_2}\loc_{v}\big(\Sel_2(J_b^{F})\big)=4.$ We conclude that  $\delta \in  \Sel_2(J_b^{F'})$ and that, by Proposition \ref{prop:variation_of_selmer_structure} (ii), we have
\[\dim_{\mathbb{F}_2} \Sel_2(J_b^{F'})=\dim_{\mathbb{F}_2} \Sel_2(J_b^{F}) -4.\]

We now relabel $F'$ as $F$  and iterate the above argument, eventually arriving at a quadratic extension $F'/k$,  unramified at $w_1$ and $w_2$, such that $\delta \in  \Sel_2(J_b^{F'})$ and $\dim_{\mathbb{F}_2}\Sel_2(J_b^{F'}) \in \{1,3\}$. If $\dim_{\mathbb{F}_2}\Sel_2(J_b^{F'})=1$ then we are done, so suppose that $\dim_{\mathbb{F}_2}\Sel_2(J_b^{F'}) =3$. 
We now once more relabel $F'$ as $F$.
Let $T=\{\alpha_0=\delta,\alpha_1,\alpha_2\}$ be an $\F_2$-basis of $\Sel_2(J_b^{F})$. 
By Lemma \ref{arranging_linearly_independent} we can assume that 
$\alpha_0,\alpha_1,\alpha_2$ are $\mathbb{F}_4$-linearly independent.  

Suppose that $F'/k$ is such that $\Sel_2(J_b^{F'}) =\Sel_2(J_b^{F})\subset\H^1(k,G)$,
and $\left \langle x,y\right \rangle^{\textup{CT}}_{F',2}=B(x,y)$, where $B(x,y)$ is an 
alternating bilinear pairing given by
\[B(\delta,\alpha_1)=B(\delta,\alpha_2)=0\quad \textup{and}\quad B(\alpha_1,\alpha_2)=1.\]
In particular, the kernel of $B$ is $\{0,\delta\}$. 
Assuming finiteness of $\Sha(J_b^{F'})[2^\infty]$, this implies that $\textup{Kum}(J_{b,\delta})$ has a $k$-point. Indeed, we have $\dim_{\mathbb{F}_2}\Sel_2(J_b^{F'})=3$ while 
$\dim_{\mathbb{F}_2}\Sha(J_b^{F'})[2]$ is even. From the exact sequence  \eqref{eq:selmer_sequence}, some element of $J_b^{F'}(k)$ must have a non-trivial image $\alpha$ in $\Sel_2(J_b^{F'})$. Then $\alpha$ is in the kernel of $\left \langle~,~\right \rangle^{\textup{CT}}_{F',2}$, hence $\alpha=\delta$ and so $\delta$ maps to $0$ in $\Sha(J_b^{F'})$.  

To construct the sought extension $F'/k$, choose cocycles $\tilde \alpha_0,\tilde \alpha_1,\tilde \alpha_2$ representing 
$\alpha_0,\alpha_1,\alpha_2$, respectively. By  Lemma \ref{is_iso_vart}, since $T$ is $\mathbb{F}_4$-linearly independent, the map $\varphi_T:\textup{Gal}(k_T/k)\to G(\bar k)^3\rtimes \textup{Gal}(P)$ is an isomorphism. 
Since $T\subset\Sel_2(J_b^{F})$, the cup products 
$\alpha_i\cup_{e_2} \alpha_j$
are trivial in $\H^2(k,\mu_2)$ for each $i,j$, where $e_2$ is the Weil pairing on $G(\bar{k})$. For each $i<j$, choose $1$-cochains $\tilde{\gamma}_{ij}:\Ga\to \F_2$ such that $d\tilde{\gamma}_{ij}=\tilde{\alpha}_i\cup_{e_2} \tilde{\alpha}_j$.  

Let $K_T\subset\bar k$ be the fixed field of $\Ker(\hat{\varphi}_T)$ as defined in 
Proposition \ref{rosenmontag}. By this proposition, we can find $\si\in\Gal(K_T/k)$
that acts on the roots of $P(t)$ as a $5$-cycle and, moreover, $\gamma_{i,j}(\si)$
has arbitrary prescribed values in $\F_2$, for all $i<j$.
Then we have $G(\bar{k})/(\sigma-1)=0$, hence for $i=0,1,2$ there is a unique element $P_{i,\sigma}\in G(\bar k)$ such that $\tilde \alpha_i(\sigma)=\si P_{i,\sigma}-P_{i,\si}$. 
We choose $\si$ so that for $i<j$ we have 
$$\tilde{\gamma}_{i,j}(\si)=e_2(P_{i,\sigma},(\sigma-1)P_{j,\sigma})+
\left \langle\alpha_i,\alpha_j\right \rangle^{\textup{CT}}_{F,2}+B(\alpha_i,\alpha_j).$$
Fix a finite set $S$ of places of $k$ containing $S_{0,b}\cup\Ram(K_TF/k)$.  Combining the last statement of Proposition \ref{rosenmontag} with Lemma \ref{is_iso_vart}, we see that any homomorphism $\textup{Gal}(K_T/k)\to \mathbb{Z}/2$ factors through $\textup{Gal}(P)$, hence vanishes on $\sigma$. Thus, arguing as in Lemma \ref{lem:existence_of_extensions}, 
we find a place $v\notin S$ such that $\Frob_v\in \textup{Gal}(K_T/k)$ is in the conjugacy class of $\sigma$, and $a\in k^\times$ such that $k(\sqrt{a})/k$ is ramified at $v$, split at all places in $S$, and is unramified elsewhere. Write $F=k(\sqrt{a_0})$ and set $F'=k(\sqrt{aa_0})$.
The $2$-Selmer conditions for $F$ and $F'$ agree at all places away from $v$. Further, at the place $v$ we have $G(\bar{k})/(\textup{Frob}_v-1)=0$, hence $\H^1(k_v,G)=0$. Thus $\Sel_2(J_b^{F'}) =\Sel_2(J_b^{F})$ as subgroups of $\H^1(k,G)$. By Proposition \ref{prop:variation_of_CTP} we have 
$\left \langle\alpha_i,\alpha_j\right \rangle^{\textup{CT}}_{F',2}=B(\alpha_i,\alpha_j)$
 for all $i<j$. Both pairings $B$ and $\left \langle~,~\right \rangle^{\textup{CT}}_{F',2}$ are alternating, so they are equal. Thus $\delta$ is the unique non-trivial element in the kernel of $\left \langle~,~\right \rangle^{\textup{CT}}_{F',2}$, so $\delta$ maps to $0$ in $\Sha(J_b^{F'})$.  
 \hfill $\Box$

Gonville and Caius College, Trinity Street, Cambridge, CB2 1TA, United Kingdom 

\texttt{ajm269@cam.ac.uk}

\bigskip

 Department of Mathematics, 
South Kensington Campus,
Imperial College London
SW7~2AZ United Kingdom \  \ and \  \ 
Institute for the Information Transmission Problems,
Russian Academy of Sciences,
Moscow 127994 Russia

\texttt{a.skorobogatov@imperial.ac.uk}

\end{document}